\documentclass[11pt,reqno,a4paper]{amsart}
\newtheorem{thm}{Theorem}[section]
\newtheorem{defi}{Definition}[section]

\newtheorem{cor}{Corollary}[section]
\newtheorem{pr}{Proposition}[section]
\theoremstyle{definition}
\newtheorem{rem}{Remark}[section]

\newcommand{\be}{\begin{equation}}
\newcommand{\ee}{\end{equation}}
\newcommand{\bea}{\begin{eqnarray}}
\newcommand{\eea}{\end{eqnarray}}
\newcommand{\beb}{\begin{eqnarray*}}
\newcommand{\eeb}{\end{eqnarray*}}
\usepackage{amssymb,amsfonts,amsthm,setspace,indentfirst,mathrsfs}
\usepackage{xcolor}
\usepackage{minibox}

\usepackage{enumitem}
\usepackage{cite}
\numberwithin{equation}{section}
\begin{document}
%
\title[Symmetry and Pseudosymmetry properties with Ricci Soliton of the RNdS spacetime]{\bf{Symmetry and Pseudosymmetry properties with Ricci soliton of the Reissner-Nordstr\"{o}m-de Sitter spacetime }}


\author[A. A. Shaikh and Kamiruzzaman]{Absos Ali Shaikh$^{* 1}$ and Kamiruzzaman$^{2}$}
\date{\today}

\address{\noindent$^{1,2}$ Department of Mathematics,
	\newline The University of Burdwan, Golapbag,
	\newline Burdwan-713104, West Bengal, India}
\email{aashaikh@math.buruniv.ac.in$^{1}$ ; aask2003@yahoo.co.in$^{1}$ }
\email{kamiruzzaman8145@gmail.com$^{2}$}


%
%
\dedicatory{}

\begin{abstract}	
	The primary objective of the article is to investigate the symmetry and pseudosymmetry properties of the Reissner-Nordström-de Sitter (briefly, RNdS) spacetime. The secondary aim of the paper is to explore the notion of Ricci solitons in RNdS spacetimes. The study is important due to the conceding of almost Ricci soliton and almost  Ricci Yamabe soliton of the RNdS spacetime. The analysis shows that this spacetime satisfies multiple types of symmetric and pseudosymmetric conditions. It is interesting to note that RNdS spacetime reveled pseudosymmetry, Weyl conformal pseudosymmetry, Weyl projective pseudosymmetry, conharmonic pseudosymmetry, and concircular pseudosymmetry. Also the RNdS spacetime is pseudosymmetric due to conformal, projective, conharmonic curvature tensors. Furthermore, in the RNdS spacetime, obtained by second order covariant derivatives $R\cdot R$ is linearly dependent on $Q(S,R)$ and $Q(g,C)$. It is demonstrated that the RNdS spacetime is 2-quasi Einstein and an Ein(2) space with recurrent conformal 2-forms. We derive the general form of the compatible tensors for this spacetime. The energy-momentum tensor of the RNdS spacetime is also shown to be pseudosymmetric and also the energy momentum tensor is pseudosymmetric due to conformal, conharmonic, concircular and projective curvature tensor. The energy momentum tensor is compatible with projective, conharmonic, concircular, Riemann and conformal curvature. We study a generalized notion of curvature inheritance and find that, with respect to the non-Killing vector fields $\partial/\partial r$ and $\partial/\partial \theta$, the RNdS spacetime does not satisfy these inheritance conditions. However, the RNdS spacetime is shown to admit an almost Ricci soliton and an almost $\eta$-Ricci Yamabe soliton with respect to the non-Killing vector field $\partial/\partial r$ but with respect to the non-Killing vector field $\partial/\partial \theta$ the spacetime does not admit such notions. Finally, we present a comparison between the RNdS and Vaidya-Bonner-de Sitter (VBdS) spacetimes.

\end{abstract}

%

\noindent\footnotetext{ $^*$ Corresponding author(Absos Ali Shaikh, E-mail address: aashaikh@math.buruniv.ac.in, aask2003@yahoo.co.in).\\
$\mathbf{2020}$\hspace{5pt}Mathematics\; Subject\; Classification: 53B20; 53B30; 53B50; 53C15; 53C25; 53C35; 83C15. }
\maketitle
%

\section{\bf Introduction}\label{intro}

Let $(\mathcal{V}, g_{ab})$ be an $n$-dimensional $(n \ge 3)$ manifold which is an open, connected  subset of $\mathbb{R}^n$  endowed with a semi-Riemannian metric $g_{ab}$ of signature $(\chi,n)$ or $(n,\chi)$, where the index $\chi$ satisfies $0 \le \chi \le n-1$. When $\chi = 1$ the metric is Lorentzian, while $\chi = 0$ represents that the metric is Riemannian. A spacetime will mean a connected $4$-dimensional Lorentzian manifold. We denote the Riemann curvature tensor by $R_{hijk}$, the Ricci tensor by $Ric_{ab}$, the scalar curvature by $\kappa$, and the Levi-Civita connection by $\nabla$, here the indices $i, j,k, h, a, b $ belongs to $\{1,2, 3,..n\}$. In addition, we use $C$, $P$, $W$ and $K$ for the conformal, projective, concircular and conharmonic curvature tensors respectively. Unless explicitly stated otherwise, the symbol $\mathcal{V}$ $(n \ge 3)$ will always refer to as a semi-Riemannian manifold.

In the framework of general relativity, the energy-momentum tensor of type $(0, 2)$ plays a central role in characterizing the physical sources within a given spacetime. This tensor significantly impacts the geometric structure by influencing quantities such as the Ricci tensor, the scalar curvature, and the metric itself. These elements collectively determine the curvature and the overall configuration of the manifold, which are essential for analyzing its geometric behavior of the manifold. Notably, the scalar curvature is obtained as the trace of the Ricci tensor, which in turn is derived as the trace of the Riemann curvature tensor. As such, the concept of “curvature” encapsulates key physical features of spacetime. For example, a Brinkmann wave qualifies as a pp-wave if it satisfies the condition $R^{ab}_{\ \ hi}R_{abjk} = 0$ \cite{Brink1925, SBH21, MS2016, SG1986}. Furthermore, spacetimes with Weyl tensors that satisfy pseudosymmetric conditions fall under Petrov type D \cite{DHV2004}. It has been shown in \cite{C01, SKH11, SYH09} that a quasi-Einstein manifold can be interpreted as a perfect fluid solution and vice versa. A notable example is the Friedmann-Lemaître-Robertson-Walker (FLRW) spacetime \cite{DDHKS00}, which exemplifies a quasi-Einstein geometry and serves as a fundamental cosmological model due to its assumptions of homogeneity and isotropy.

In contrast, the geometry of a semi-Riemannian manifold is primarily governed by its curvature, which is intricately related to the covariant derivatives of various tensors, possibly of higher order. When the condition $\nabla_r R_{hijk} = 0$ holds, the manifold $\mathcal{V}$ is said to be locally symmetric \cite{Cart26}, indicating that the geometry around each point is symmetric, and the corresponding local geodesic symmetries are isometries. This concept forms a cornerstone of differential geometry, emphasizing the consistent geometric behavior in the neighborhood of each point. Ruse \cite{Patt52, Ruse46, Ruse49a, Ruse49b} extended this idea by introducing the notion of “kappa” spaces, which generalize local symmetry by relaxing the strict condition on curvature, thereby encompassing a wider range of geometric structures. Initially developed through the study of curvature preservation under parallel transport, these “kappa” spaces permit scalar scaling and have since been acknowledged as an important class of manifolds. Walker \cite{Walk50} later termed such structures as recurrent manifolds, leading to further developments and the classification of numerous generalized recurrent types of manifold. For readers seeking an in depth understanding of these families of manifolds, references such as \cite{SAR13, SK14, SKA18, SP10, SR11, SRK15, SRK17, SRK18} offer extensive analyses. Earlier, Shirokov \cite{Shi25} explored manifolds characterized by the presence of covariantly constant tensors, with additional elaboration provided in the subsequent work \cite{Shi98}.

The integrability conditions arising from the equations $\nabla R = 0$ are succinctly expressed by the formula $R \cdot R = 0$. These conditions have been extensively studied in the works of Cartan \cite{Cart46} and Shirokov \cite{Shi25}, where they are analyzed in depth. When $R \cdot R = 0$, the manifold $\mathcal{V}$ is named as semisymmetric. The specific conditions under which $R \cdot R = L_{R} Q(g, R)$ are thoroughly explored in Sinyukov’s work \cite{Sin54}, where $L_{R}$ denotes any smooth function defined on the subset 
$\{x \in M : Q(g, R) x \neq 0\}$, with $Q(g, R)$ representing the Tachibana tensor. These conditions are referred to as generalized symmetric and are explained in detail in \cite{MIKES76} and \cite{Add1}, with more concise equations provided in \cite{MIKES88, MIKES92, MIKES96, MSV15, MVH09}.

The study of symmetries in differential geometry and general relativity has evolved significantly over the past few decades, with important contributions enhancing our understanding of geometric structures in curved spacetimes. A notable advancement was made by Adamów and Deszcz \cite{AD83}, who introduced the concept of second-order symmetry, extending the classical notion of semisymmetry. Their work, which emerged from a detailed investigation of Einstein’s field equations and the geometry of totally umbilical hypersurfaces, laid the groundwork for further exploration into generalized symmetry conditions.

One of the key outcomes of this line of research is the notion of pseudosymmetry, which has since found relevance in a wide range of spacetimes. For instance, the Robertson-Walker cosmological models \cite{ADEHM14, DK99}, the Schwarzschild solution \cite{Sch1916}, and the Reissner-Nordström black hole geometry \cite{Reis1916, Nord1918} are all recognized as classical examples of pseudosymmetric manifolds. The geometric implications and distinctions among local symmetry, semisymmetry, and pseudosymmetry are comprehensively examined by Haesen and Verstraelen in a series of works \cite{HV07, HV07A, HV09}, which serve as valuable references for understanding the theoretical underpinnings of these conditions.

Deszcz's formulation of pseudosymmetry has remained a central focus in both differential geometry and general relativity for the last four decades. It has been widely applied in general relativity and cosmology to characterize various spacetime models. An extensive list of spacetimes exhibiting pseudosymmetric curvature structures spanning different physical interpretations and geometric characteristics can be found in studies such as \cite{ADEHM14, DK99, Kowa06, SAA18, SAAC20, SDKC19, SAACD_LTB_2022, SADM_TCEBW_2023, EDS_sultana_2022, ESKiselev2024, SADM_TCEBW_2023, SAS_EiBI_2024, SASK_pgm_2024, SHDS_hayward, SHS_warped, SHS_Bardeen, SHSB_VBDS}.

Complementing this research, Mikeš and his collaborators \cite{MIKES76, MIKES88, MIKES92, MIKES96, MS94, MSV15, MVH09} conducted a thorough analysis of geodesic mappings in symmetric Riemannian manifolds, offering further insight into curvature-preserving transformations and their role in the classification of manifold structures. Furthermore, the relationship between Deszcz pseudosymmetry and Chaki pseudosymmetry has been critically examined in \cite{SDHJK15}, highlighting areas of overlap and contributing to a more nuanced understanding of symmetry concepts in pseudo-Riemannian geometry.
Together, these works provide a robust foundation for current investigations into curvature-restricted geometries and continue to inform contemporary approaches to modeling spacetime in mathematical physics. \\

Black holes are fascinating cosmic entities. The interior of the black hole is bounded by an event horizon, beyond which neither light nor matter can return once it has fallen in, and hence black holes can only absorb material. Yet, due to quantum effects, they are capable of emitting Hawking radiation and thereby losing energy \cite{Hwkn1975}.
Black holes that possesses electric charge display distinct characteristics, as the presence of charge notably influences both their geometric structure and thermodynamic behavior.
The discovery of black hole solutions to Einstein’s field equations in the early 1960s, along with recent breakthroughs such as the observation of gravitational waves and the first direct image of a black hole’s event horizon in the nearby galaxy M87, represent significant milestones that carry deep philosophical significance.
Black holes are gravitationally collapsed entities that exhibit infinite redshift, rendering them causally isolated from the rest of the universe. Any event taking place within a black hole has no effect on the outside world. Therefore, one of the fundamental challenges is to clearly define the boundary separating the interior of a black hole from its exterior.
Since their discovery, black holes have inspired scientists to explore gravity more deeply. While investigating their physical characteristics is still difficult, the theoretical framework of black holes has raised important questions and led to significant predictions about space, time, and the cosmos. These explorations have contributed to the development of theories aiming to unify the fundamental forces of nature. Because any viable theory of quantum gravity must reduce to established physics under specific conditions and since black holes lie at the crossroads of classical and quantum phenomena, the black holes described by such theories are especially noteworthy.

To begin, it is essential to examine the properties of black holes as described by general relativity (GR), which is our prevailing theory of gravity. Most black hole studies are based on the assumption of a universe with three spatial and one temporal dimension—commonly referred to as 3+1 dimensional spacetime—consistent with current observational evidence. However, many unification theories suggest the existence of more than four dimensions, making it valuable to explore and classify higher-dimensional black hole solutions within GR for comparison. At a first glance, one might expect black holes in higher dimensions to be straightforward extensions of their 3+1 dimensional counterparts, which justifies a brief review of black holes in four-dimensional spacetime.

A fundamental aspect in the classification of black holes is the principle of uniqueness. Research (e.g., \cite{bb1}) has demonstrated that, under fairly general assumptions like asymptotic flatness and stationarity, black holes in $3+1$ dimensions are fully characterized by just three parameters: mass, angular momentum, and electric charge. The Kerr–Newman (KN) solution \cite{bb2} to the Einstein–Maxwell equations in four-dimensional spacetime represents the most comprehensive black hole model, as it incorporates all three of these quantities. Earlier solutions namely, the Schwarzschild solution \cite{Sch1916} (describing mass alone), the Reissner–Nordström solution \cite{Reis1916, bb5} (mass and charge), and the Kerr solution \cite{bb6} (mass and angular momentum) can be seen as special or limiting cases of the KN black hole. As a result, analyzing the KN solution provides insight applicable to all four-dimensional black holes. Importantly, these black hole geometries feature a null hypersurface, topologically equivalent to a two-sphere ($S^2$) \cite{bb7}, known as the event horizon. This surface functions as a one-way boundary that separates the black hole's interior from the external universe. The interaction between a black hole’s mass and the surrounding anisotropic fluid is a phenomenon that has also been explored in \cite{CI2025}.

The stationary, spherically symmetric\cite{HAli2003} solutions to the Einstein-Maxwell gravity with positive cosmological constant\cite{GG2019} is known as Reissner-Nordstr\"{o}m-de Sitter spacetime [\cite{CI2025}, Section 2, equation 2--4].

The metric of a Reissner-Nordstr\"{o}m-de Sitter spacetime is as follows
\begin{eqnarray}\label{RNdSM}
ds^2= -f(r)dt^2+\frac{dr^2}{f(r)}+r^2( d\theta^2+\sin^2\theta d\phi^2),
\end{eqnarray}
where  $$f(r)=1-\frac{2M}{r}+\frac{e^2}{r^2}-\Lambda r^2,$$ with $e$ being the electric charge (a constant), $M$ is the Mass of the black hole, $\Lambda$ is the cosmological constant, -$\infty < t < \infty,$ $r \geq 0$ and $\phi \in  [0, 2 \pi)$  defining the temporal, radial, and angular coordinates, respectively.

In particular\\
$(a)$ if $e \neq 0, \Lambda=0$ the metric \eqref{RNdSM} becomes Reissner-Nordstr\"{o}m spacetime,\\
$(b)$ if $e=0 , \Lambda \neq 0$ the metric \eqref{RNdSM} turns into Kottler spacetime,\\
$(c)$ if $e=0, \Lambda=0$ the metric \eqref{RNdSM} reduces to Schwarzschild spacetime respectively \cite{Kowa06}.
These spacetimes serve as examples of non-semisymmetric and pseudosymmetric manifolds, as illustrated in [\cite{DVV1991}, Example 1]. It is well established that the Kottler spacetime represents a non-Ricci-flat Einstein manifold, whereas the Schwarzschild spacetime is characterized by being Ricci-flat.

During last few decades, several physical phenomenons, such as charged particle and strong cosmic censorship \cite{GG2019}, spinning particles \cite{HAli2003}, Hawking–Rényi black hole thermodynamics, Kiselev solution,
and cosmic censorship \cite{CI2025},  bounds with
quasinormal frequencies \cite{CCDNP2023}, critical phenomena and thermodynamics \cite{Zhao2014}  and Schottky anomaly \cite{Zhen2025} have been investigated on the RNdS spacetime by various physicists.

 In the existing literature, Kowalczyk\cite{Kowa06} investigated certain properties of Reissner-Nordstr\"{o}m-de Sitter type spacetimes. In the present article, we examine several intriguing curvature and Ricci soliton‑related features that have not been discussed previously in the literature.

While many researchers have explored curvature-restricted geometric structures in various spacetimes, the RNdS spacetime remains relatively underexamined in this context. The central aim of this study is to analyze the curvature characteristics of the RNdS spacetime in detail. Our investigation reveal that this spacetime exhibits multiple forms of pseudosymmetry, including those associated with the conformal, conharmonic, concircular, Weyl projective curvature tensors, as well as Ricci pseudosymmetry. Notably, the RNdS spacetime does not satisfy any semisymmetric conditions.

Furthermore, it is identified as a Roter-type manifold, an Einstein manifold of level $2,$ and a $2$-quasi-Einstein manifold. It is also observed that this spacetime admits neither curvature collineations nor curvature inheritance (including generalized curvature inheritance).
The most interesting observation about RNdS spacetime is that it does not admit any of the geometric structure obtained by the restriction of first order covariant derivatives, such as local symmetry, recurrency, conformal symmetry etc. But RNdS spacetime obeys all the curvature properties produced by second order covariant derivatives.\\

The  collection of all Killing vector fields on the manifold $\mathcal{V}$, denoted by $\mathcal{K}(\mathcal{V})$, constitutes a Lie subalgebra of $\chi(\mathcal{V})$. It is established that $\mathcal{K}(\mathcal{V})$ can contain at most $\frac{n(n+1)}{2}$ linearly independent Killing vector fields. If this maximum number is achieved, meaning $\mathcal{K}(\mathcal{V})$ has exactly $\frac{n(n+1)}{2}$ such vector fields, then the manifold $\mathcal{V}$ is considered to as a maximally symmetric space. A space of this kind possesses constant scalar curvature. Section $3$ demonstrates that the scalar curvature $\kappa$ of the RNdS spacetime is $\kappa = -12 \Lambda$, which remains constant (since $\Lambda$ is a cosmological constant). Therefore, the RNdS spacetime qualifies as a maximally symmetric space.

In 1982, Richard Hamilton introduced the idea of Ricci flow as a tool to investigate the deformation of Riemannian metrics over time, specifically on compact three-dimensional manifolds with positive Ricci curvature \cite{Hamilton1982}. A special class of solutions to the Ricci flow that evolve in a self-similar solution are called Ricci solitons. These geometric structures generalize the notion of Einstein metrics and have attracted considerable attention in both  differential geometry and mathematical physics \cite{Bess87, Brink1925, SYH09, S09}. The original concept has since been broadened to include several variants such as almost Ricci solitons, $\eta$-Ricci solitons, and almost $\eta$-Ricci solitons, each providing deeper insight into the geometric and topological behavior of manifolds under curvature flow.

Over the last thirty years, a substantial body of research has focused on Ricci solitons, Yamabe solitons, and their numerous generalizations. Notable contributions can be found in works such as \cite{AliAhsan2013, AliAhsan2015, Ahsan2018, SDAA_LCS_2021, SMM2022, SM23, Absos23, SMB24}, among others. This growing interest has established these topics as active and influential areas within contemporary differential geometric research. The secondary aim of this article is to investigate whether RNdS spacetime admits Ricci solitons or not. The answer of this question has given as affirmative.


 Under certain conditions, the spacetime supports both almost Ricci solitons and almost $\eta$-Ricci–Yamabe solitons associated with non-Killing vector fields. Lastly, a comparative analysis is performed between the geometric structures of the RNdS spacetime and the Vaidya–Bonner–de Sitter (VBdS) spacetime.\\

The paper is structured as follows. In Section $2$, we discuss the basic definitions of the geometric frameworks needed to study the RNdS spacetime. Section $3$ presents an in-depth analysis of the RNdS spacetime, highlighting several noteworthy findings. In Section $4,$ we investigate the geometric characteristics tied to the energy--momentum tensor of the RNdS spacetime. Section $5$ explores the existence and properties of Ricci and Ricci--Yamabe solitons within the RNdS spacetime. Finally, Section $6$ offers a comparative study of the geometric structures in the RNdS and VBdS spacetimes.
\section{\bf Curvature restricted geometric properties}
Given two symmetric $(0,2)$ tensors $\mathcal{F}$ and $\mathcal{H}$, their Kulkarni–Nomizu product, denoted by $(\mathcal{F} \wedge \mathcal{H})$, results as a $(0,4)$ tensor. This product, often used to describe specific curvature-related characteristics, is defined as follows (refer to \cite{DGHS11, DHJKS14, Glog02, G08, Kowa06}):

$$(\mathcal{F}\wedge \mathcal{H})_{ijkl}=2\mathcal{F}_{i[l}\mathcal{H}_{k]j} +2\mathcal{F}_{j[k}\mathcal{H}_{l]i},$$
Here, the notation $[.]$ represents antisymmetrization over the enclosed indices, which guarantees that the expression is antisymmetric with respect to those index pairs.

The expression of second kind Christoffel symbols $(\Gamma^i_{jk} )$ is defined as follows:\\
Given a coordinate basis $(x^1, x^2, \ldots, x^n)$ and a metric tensor $g_{ij}$, the Christoffel symbols of the second kind are given by:

$$
\Gamma^k_{ij} = \frac{1}{2} g^{kl} \left( \frac{\partial g_{jl}}{\partial x^i} + \frac{\partial g_{il}}{\partial x^j} - \frac{\partial g_{ij}}{\partial x^l} \right).
$$

The $(1,3)$-type curvature tensors namely the Riemann, conharmonic, concircular, conformal, and projective curvatures are each defined in the following manner:
\bea
R^h_{ijk}&=& 2\left(\Gamma^s_{i[j}\Gamma^h_{k]s} + \partial_{[k}\Gamma^h_{j]i}\right),\nonumber\\
K^h_{ijk}&=& R^h_{ijk} - \frac{2}{n-2}\left( Ric^h_{[i}g_{j]k} + \delta^h_{[i}Ric_{j]k}\right),\nonumber \\
W^h_{ijk}&=& R^h_{ijk} - \kappa\frac{2}{n(n-1)}\delta^h_{[i}g_{j]k}, \nonumber \\
C^h_{ijk}&=& R^h_{ijk}+\frac{2}{n-2}\left( \delta^h_{[i}Ric_{j]k} +Ric^h_{[i}g_{j]k}\right) -\kappa\frac{2}{(n-1)(n-2)}\delta^h_{[i}g_{j]k} ,\nonumber \\
P^h_{ijk}&=& R^h_{ijk} -\frac{2}{n-1}\delta^h_{[i}Ric_{j]k}, \nonumber
\eea
here, $\Gamma^i_{jk}$ stands for the connection coefficients, $Ric^i_j$ indicates the Ricci tensor of type $(1,1)$, and $\partial_h$ denotes the partial derivative operator $\frac{\partial}{\partial x^h}$.

The tensors of rank $(0,4)$ namely $R_{hijk},$ $K_{hijk},$ $W_{hijk},$ $C_{hijk},$ and $P_{hijk}$ are derived by lowering indices with the help of the metric tensor $g_{hi}$. More precisely, they are given by:

\bea
R_{hijk}&=& g_{h\alpha}(\partial_k \Gamma^\alpha_{ij}-\partial_j \Gamma^\alpha_{ik}+\Gamma^\beta_{ij}\Gamma^\alpha_{\beta k}-\Gamma^\beta_{ik}\Gamma^\alpha_{\beta j}) , \nonumber \\
K_{hijk}&=& R_{hijk} - \frac{1}{n-2} (g\wedge Ric)_{hijk} , \nonumber \\
W_{hijk}&=& R_{hijk} - \frac{\kappa}{2n(n-1)} (g\wedge g)_{hijk}, \nonumber \\
C_{hijk}&=& R_{hijk}-\frac{1}{n-2}(g\wedge Ric)_{hijk}+\frac{\kappa}{2(n-1)(n-2)}(g\wedge g)_{hijk} ,\nonumber \\
P_{hijk}&=& R_{hijk} -\frac{1}{n-1}(g_{hk}Ric_{ij}-g_{ik}Ric_{hj}). \nonumber
\eea

Let $Y$ be a tensor of type $(0,d)$ with $d \geq 1$. The tensors of type $(0,d+2)$, denoted by $(X \cdot Y)$ (refer to \cite{DG02, DGHS98, DH03, SDHJK15, SK14}) and $Q(\mathscr{Y}, Y)$ (see \cite{DGPSS11, SDHJK15, SK14, Tach74}), are defined in the following way:

\beb
(X\cdot Y)_{b_1b_2\cdots b_dws}&=&-\left[ X^\alpha_{ksb_1}Y_{\alpha b_1\cdots b_d}+ \cdots + X^\alpha_{ksb_d}Y_{b_1\cdots \alpha}\right], \\
Q(\mathscr{Y},Y)_{b_1b_2\cdots b_d ij}&=&\mathscr{Y}_{jb_1}Y_{ib_2\cdots b_d}+ \cdots + \mathscr{Y}_{ib_d}Y_{b_1b_2\cdots j}\\ 
&-& \mathscr{Y}_{ib_1}B_{jb_2\cdots b_d}- \cdots - \mathscr{Y}_{ib_d}Y_{b_1b_2\cdots j},
\eeb
where $\mathscr{Y}_{ij}$ is the $(0,2)$ rank symmetric tensor and $X^h_{ijk}$ is the $(1,3)$ rank tensor.\\
\begin{defi} \cite{AD83, Cart46, Chak87, Chak88, Desz92, Desz93, DGHZ15, DGHZ16,SAAC20N, SK14, SKppsnw, Szab82, Szab84, Szab85} 

A manifold $\mathcal{V}$ is called $Y$-pseudosymmetric with respect to the tensor $X$ if the expression $X \cdot Y$ is proportional to $Q(g, Y)$, that is, $X \cdot Y = f_Y Q(g, Y)$, where $f_Y$ is a smooth scalar function defined on $\mathcal{V}$. Moreover, if the relation $X \cdot Y = \bar{f}_Y Q(S, Y)$ holds, where $\bar{f}_Y$ is also a smooth function on $\mathcal{V}$, then the manifold is said to be Ricci-generalized $Y$-pseudosymmetric with respect to $X$.
\end{defi} 

When both $X$ and $Y$ are chosen as the Riemann curvature tensor $R$, the manifold is described as pseudosymmetric or Ricci-generalized pseudosymmetric. By examining different curvature tensors such as the Riemann, Ricci, and others various forms of pseudosymmetric and Ricci-generalized pseudosymmetric manifolds can be classified.

A manifold is called $Y$-semisymmetric with respect to a tensor $X$ if it satisfies the condition $X \cdot Y = 0$. In particular, a semisymmetric manifold is one for which $R \cdot R = 0$ holds \cite{Cart46, Szab82, Szab84, Szab85}. Semisymmetric manifolds form a special class within the broader category of pseudosymmetric manifolds, emphasizing the significance of specific curvature relations in differential geometry. This classification illustrates how conditions like $X \cdot Y = 0$ help to define distinct symmetric and pseudosymmetric properties in the geometric structure of manifolds.

\begin{defi}$($\cite{C01,DGHZ16, DGJZ-2016, DGP-TV-2011, S09, SKH11, SK19,SYH09}$)$ 
 
A manifold $\mathcal{V}$ is called a $d$-quasi-Einstein manifold if, for some scalar $L$ and $0 \leq d \leq (n-1)$, the tensor $(Ric - L g)$ has rank $d$. In particular, when $d = 1$ (respectively, $d = 0$), the manifold is known as quasi-Einstein \cite{SKH11, SYH09} (respectively, Einstein). \\
Additionally, if $L = 0$ in the quasi-Einstein condition, the manifold is referred to as Ricci simple.
 
\end{defi}

It is worth mentioning that the Robertson–Walker spacetime \cite{ARS95, Neill83, SKMHH03} falls under the category of quasi-Einstein manifolds, while the Kaigorodov spacetime \cite{SDKC19} is an example of an Einstein manifold. The Kantowski–Sachs spacetime \cite{SC21} and the Som–Raychaudhuri spacetime \cite{SK16srs} serve as instances of 2-quasi-Einstein manifolds. Furthermore, spacetimes such as Vaidya \cite{SKS19}, Gödel \cite{DHJKS14}, and Morris–Thorne \cite{ECS22} are classified as Ricci simple manifolds.

\begin{defi}
A manifold $\mathcal{V}$ is said to have a cyclic parallel Ricci tensor (see \cite{Gray78, SB08, SJ06, SJ07}) if the following condition is satisfied throughout $\mathcal{V}$:
\[ (\nabla_{\mathcal{D}_1} Ric)(\mathcal{D}_2, \mathcal{D}_3) + (\nabla_{\mathcal{D}_2} Ric)(\mathcal{D}_3, \mathcal{D}_1) + (\nabla_{\mathcal{D}_3} Ric)(\mathcal{D}_1, \mathcal{D}_2) = 0 \]
where $\mathcal{D}_1, \mathcal{D}_2, \mathcal{D}_3$ are arbitrary vector fields. Moreover, the Ricci tensor is known as a Codazzi type (see \cite{F81, S81}) if it meets the criterion
\[ (\nabla_{\mathcal{D}_1} Ric)(\mathcal{D}_2, \mathcal{D}_3) = (\nabla_{\mathcal{D}_2} Ric)(\mathcal{D}_1, \mathcal{D}_3) \]
for all vector fields $\mathcal{D}_1$, $\mathcal{D}_2$, and $\mathcal{D}_3$ on $\mathcal{V}$.
\end{defi}
It is important to highlight that the Ricci tensor in the $(t\text{-}z)$-plane wave spacetime \cite{EC21} is of Codazzi type, whereas cyclic parallel Ricci tensor behavior has been observed in the Gödel spacetime \cite{DHJKS14}.

\begin{defi} $($\cite{Bess87, SK14, SK19}$)$ 
A manifold $\mathcal{V}$ is referred to as an Einstein manifold of level $4$ if it fulfills the condition

$$
\bar{\gamma}_1 Ric^4 + \bar{\gamma}_2 Ric^3 + \bar{\gamma}_3 Ric^2 + \bar{\gamma}_4 Ric + \bar{\gamma}_5 g = 0, \quad (\bar{\gamma}_1 \neq 0),
$$

where each $\bar{\gamma}_i$ is a smooth function defined on $\mathcal{V}$. If $\bar{\gamma}_1 = 0$ but $\bar{\gamma}_2 \neq 0$, the manifold is classified as an Einstein manifold of level $3$, similarly, if both $\bar{\gamma}_1 = \bar{\gamma}_2 = 0$ and $\bar{\gamma}_3 \neq 0$, then it is considered an Einstein manifold of level $2$.

\end{defi}

It has been observed that the Vaidya–Bonner spacetime \cite{SDC} and the Lifshitz spacetime \cite{SSC19} serve as examples of level 3 Einstein manifolds, while the Siklos spacetime \cite{SDKC19}, Lemos black hole spacetime \cite{SASK_LBH_2025} and the Nariai spacetime \cite{SAAC20N} represent level 2 Einstein manifolds.

\begin{defi} 

A manifold $\mathcal{V}$ is called generalized Roter type if its Riemann curvature is expressed as \cite{Desz03, DGJPZ13, DGJZ-2016, DGP-TV-2015, SK16, SK19}:

$$
R = (\mathcal{L}_{11} Ric^2 + \mathcal{L}_{12} Ric + \mathcal{L}_{13} g) \wedge Ric^2 + (\mathcal{L}_{22} Ric + \mathcal{L}_{23} g) \wedge Ric + \mathcal{L}_{33} (g \wedge g),
$$

where $\mathcal{L}_{ij}$ are scalar functions. If these coefficients are such that $R$ can be represented as a linear combination involving only $g \wedge g$, $g \wedge Ric$, and $Ric \wedge Ric$, then the manifold is referred to as a Roter type manifold \cite{Desz03, DG02, DGP-TV-2011, DPSch-2013, Glog-2007}.

\end{defi}

It is important to note that the Nariai spacetime \cite{SAAC20N}, Melvin magnetic spacetime \cite{SAAC20}, and Robinson–Trautman spacetime \cite{SAA18} are identified as Roter type manifolds. In contrast, the Lifshitz spacetime \cite{SSC19} and Vaidya–Bonner spacetime \cite{SDC} are examples of generalized Roter type manifolds.

%
%
%

\begin{defi}$($\cite{DD91,DGJPZ13, MM12a, MM12b, MM13,MM22b}$)$
Let $\mathcal{B}_{\varpi}$ denote the endomorphism corresponding to a symmetric $(0,2)$ tensor $\varpi$, and let $B$ be a $(0,4)$ tensor field. The tensor $\varpi$ is said to be $B$-compatible if the condition

$$
\mathop{\mathcal{S}}_{\mathcal{Q}_1,\mathcal{Q}_2,\mathcal{Q}_3} B(\mathcal{B}\mathcal{Q}_1,\mathcal{Q},\mathcal{Q}_2,\mathcal{Q}_3) = 0
$$

is satisfied on the manifold $\mathcal{V}$. Moreover, if the tensor $\Psi \otimes \Psi$ satisfies this compatibility condition with respect to $B$, then the associated $1$-form $\Psi$ is called $B$-compatible.

\end{defi}

To define the corresponding compatibility conditions for curvature tensors such as $R$ (Riemann), $C$ (conformal), $W$ (concircular), $P$ (projective), and $K$ (conharmonic), one substitutes $B$ with the respective tensor in the general compatibility framework.

\begin{defi}\label{defi2.8} 
Let $B$ be a $(0,4)$ tensor field. The associated curvature $2$-forms $\Pi_{(B)}^m l$ \cite{LR89} are said to be recurrent \cite{MS12a, MS13a, MS14} if and only if the following condition holds:

$$
\mathop{\mathcal{S}}_{\mathcal{Q}_1,\mathcal{Q}_2,\mathcal{Q}_3}(\nabla_{\mathcal{Q}_1}B)(\mathcal{Q}_2,\mathcal{Q}_3,\mathcal{Q}_4,\mathcal{Q}_5) = \mathop{\mathcal{S}}_{\mathcal{Q}_1,\mathcal{Q}_2,\mathcal{Q}_3} \sigma(\mathcal{Q}_1) B(\mathcal{Q}_2,\mathcal{Q}_3,\mathcal{Q}_4,\mathcal{Q}_5),
$$

where $\sigma$ is a $1$-form.

In addition, $1$-forms $\lambda_{(T)l}$ are called recurrent \cite{SKP03} if they satisfy the relation

$$
(\nabla_{\mathcal{Q}_1}T)(\mathcal{Q}_2,\mathcal{Q}) - (\nabla_{\mathcal{Q}_2}T)(\mathcal{Q}_1,\mathcal{Q}) = \Gamma(\mathcal{Q}_1) T(\mathcal{Q}_2,\mathcal{Q}) - \Gamma(\mathcal{Q}_2) T(\mathcal{Q}_1,\mathcal{Q}),
$$

on the manifold $\mathcal{V}$, where $\Gamma$ is a $1$-form.

\end{defi}

\begin{defi}$($\cite{P95, Venz85}$)$ Consider the expression

$$
\mathop{\mathcal{S}}_{\mathcal{Q}_1,\mathcal{Q}_2,\mathcal{Q}_3} \Gamma(\mathcal{Q}_1) \otimes A(\mathcal{Q}_2, \mathcal{Q}_3, \mathcal{Q}_4, \mathcal{Q}_5) = 0,
$$

where $\mathcal{S}$ represents cyclic summation over the vector fields $\mathcal{Q}_1$, $\mathcal{Q}_2$, and $\mathcal{Q}_3$.

The collection of all $1$-forms $\Gamma$ satisfying this relation constitutes a vector space, denoted by $L(\mathcal{V})$. According to Venzi's classification, the manifold $\mathcal{V}$ is called an $A$-space if the dimension of $L(\mathcal{V})$ is at least one.
	
	\end{defi}
%

\begin{defi}If $\xi$ is a vector field on $\mathcal{V}$ such that $\pounds_\xi g = 0,$ then $\xi$ is called a Killing vector field. where $\pounds_\xi$ denotes the Lie derivative in the direction of $\xi$ .

\end{defi}

\begin{defi}\label{def_CI}
The manifold $(\mathcal{V}, \xi, g)$ is called Ricci soliton if 

$$
\frac{1}{2} \pounds_\xi g + Ric - \mu g = 0,
$$

where $\mu$ is a real constant and $\xi$ is called soliton vector field.
\end{defi}
 The classification of the Ricci soliton is determined by the sign of $\mu$, it is called expanding when $\mu < 0$, steady when $\mu = 0$, and shrinking when $\mu > 0$. If the constant $\mu$ is instead replaced by a smooth function, then the manifold $\mathcal{V}$ is referred to as an almost Ricci soliton \cite{Pigola2011}.
If the soliton vector field $\xi$ associated with a Ricci soliton is a Killing vector field, the Ricci soliton reduces to an Einstein manifold.
\begin{defi}\label{def_CI} (\cite{Cho2009})
The manifold $\mathcal{V}$ is called an $\eta$-Ricci soliton if there exists a non-zero $1$-form $\eta$, along with constants $\mu$ and $\alpha$, such that the following condition holds:

$$
\frac{1}{2} \pounds_\xi g + Ric - \mu g + \alpha (\eta \otimes \eta) = 0.
$$
\end{defi}
If the constants $\mu$ and $\alpha$ are instead replaced by smooth functions on $\mathcal{V}$, then the manifold is said to admit an almost $\eta$-Ricci soliton structure \cite{Blaga2016}.

In 1988, Hamilton \cite{Hamilton1988} introduced the Yamabe flow as a counterpart to the Ricci flow. Building upon these ideas, Güler and Crășmăreanu \cite{Guler2019} later proposed the Ricci-Yamabe flow, which unifies the characteristics of both Ricci and Yamabe flows into a single geometric evolution process. The self-similar solutions arising from this combined flow are known as Ricci-Yamabe solitons, while those specific to the Yamabe flow are referred to as Yamabe solitons.
\begin{defi}\label{def_CI} (\cite{Siddiqi2020})
The manifold $(\mathcal{V}, \xi, g)$ is said to be Ricci-Yamabe soliton if it fulfills the condition\cite{Siddiqi2020}

$$
\frac{1}{2} \pounds_\xi g + \phi_1 Ric + \left(\mu - \frac{1}{2} \phi_2 \kappa \right) g = 0,
$$

where $\phi_1$, $\phi_2$, and $\mu$ are constants and $\kappa$ represents the scalar curvature of the manifold.
\end{defi}
The introduction of Ricci-Yamabe solitons enhances the study of geometric flows by creating a unified structure that merges the features of both Ricci and Yamabe solitons. This framework provides a richer understanding of how manifolds evolve under the influence of these combined flows. The formulation of the Ricci-Yamabe flow marks a notable progression in the field of geometric analysis, extending Hamilton’s foundational contributions to new areas of exploration, particularly within the setting of semi-Riemannian geometry.

In semi-Riemannian manifolds, various choices of the constants $\phi_1$, $\phi_2$, and $\mu$ give rise to distinct types of soliton structures. For example, the condition $\phi_1 = 0$ and $\phi_2 = 2$ characterizes a Yamabe soliton, while $\phi_1 = 1$ and $\phi_2 = 0$ corresponds to a Ricci soliton. These solutions serve as fundamental models for understanding the manifold’s behavior under their respective geometric flows. Furthermore, when the parameters $\phi_1$, $\phi_2$, and $\mu$ are allowed to vary as smooth, non-constant functions on the manifold, the resulting structure is referred to as an almost Ricci-Yamabe soliton \cite{Siddiqi2020}. This extension introduces greater flexibility and enables the study of more diverse and complex geometric configurations.
\begin{defi}\label{def_CI} (\cite{Siddiqi2020})
The manifold $\mathcal{V}$ is called an $\eta$-Ricci-Yamabe soliton if there exists a non-zero $1$-form $\eta$, along with constants $\phi_1$ $\phi_2, and $ $\mu$ such that the following condition holds:
$$
\frac{1}{2} \pounds_\xi g + \phi_1 Ric + \left(\mu - \frac{1}{2} \phi_2 \kappa \right) g + \alpha_1 \, \eta \otimes \eta = 0,
$$

\end{defi}
where $\alpha_1$ is a scalar coefficient . If the constants $\phi_1$, $\phi_2$, $\mu$, and $\alpha_1$ are generalized to smooth functions on the manifold, the resulting structure is referred to as an almost $\eta$-Ricci-Yamabe soliton, which allows for the investigation of more intricate geometric features and their significance in the broader context of differential geometry.

For a $(1,3)$-type Rieman curvature tensor, the notion of curvature collineation was introduced by Katzin and collaborators in 1969 \cite{KLD1969, KLD1970} and is defined as $\pounds_\xi \widetilde{R}=0$, where $\widetilde{R}$ is $(1,3)$-type Riemann curvature tensor. Building upon this foundation, Duggal extended the idea in 1992 \cite{Duggal1992} by introducing the notion of curvature inheritance for the same type of tensor. Over the past thirty years, numerous research articles (see \cite{Ahsan1978, Ahsan1977_231, Ahsan1977_1055, Ahsan1987, Ahsan1995, Ahsan1996, Ahasan2005, AhsanAli2014, AA2012, AH1980, AliAhsan2012, SASZ2022, ShaikhDatta2022}) have explored and developed various aspects of these symmetry properties.

\begin{defi}\label{def_CI} (\cite{Duggal1992}) The manifold $\mathcal{V}$ is said to admit curvature inheritance if it satisfies the condition

$$
\pounds_\xi \widetilde{R} = \mathcal{G} \widetilde{R},
$$

where $\xi$ is a non-Killing vector field and $\mathcal{G}$ is a scalar function. The connection between the $(1,3)$-type curvature tensor $\widetilde{R}$ and the $(0,4)$-type curvature tensor $R$ is established through the relation

$$
R(\mathcal{G}_1, \mathcal{G}_2, \mathcal{G}_3, \mathcal{G}_4) = g(\widetilde{R}(\mathcal{G}_1, \mathcal{G}_2)\mathcal{G}_3, \mathcal{G}_4).
$$

In particular, when $\mathcal{G} = 0$, meaning $\pounds_\xi \widetilde{R} = 0$, the manifold satisfies the curvature collineation condition, as introduced in \cite{KLD1969, KLD1970}.




\end{defi}

 \begin{defi} \label{def_RI}(\cite{Duggal1992})
The manifold $\mathcal{V}$ is said to admit Ricci inheritance if there exists a vector field $\xi$ and a scalar function $\mathcal{G}$ such that

$$
\pounds_\xi Ric = \mathcal{G} Ric,
$$

where $Ric$ represents the Ricci tensor of $\mathcal{V}$.

Specifically, if $\mathcal{G} = 0$, the condition reduces to Ricci collineation, which means that the Lie derivative of the Ricci tensor vanishes, i.e., $\pounds_\xi Ric = 0$.

\end{defi}

%
%
%
%
%
%

\section{\bf Reissner-Nordstr\"{o}m-de Sitter spacetime admitting geometric structures} 
%
In the coordinates  $(t,r,\theta,\phi)$, the metric tensor  of RNdS spacetime is given by:
$$g=\left(\begin{array}{cccc}
	-(1+\frac{e^2}{r^2}- \frac{2M}{r}-\Lambda r^2) & 0 & 0 & 0 \\
	0 & (1+\frac{e^2}{r^2}- \frac{2M}{r}-\Lambda r^2)^{-1} & 0 & 0 \\
	0 & 0 & r^2 & 0 \\
	0 & 0 & 0 & r^2 \sin^2\theta 
\end{array}
\right).
$$
That is,
 \begin{eqnarray}
g_{11}=-\left(1+\frac{e^2}{r^2}- \frac{2M}{r}-\Lambda r^2\right),\; g_{22}=\left(1+\frac{e^2}{r^2}- \frac{2M}{r}-\Lambda r^2\right)^{-1}, \ \ g_{33}=r^2, \;g_{44}=r^2 \sin^2\theta , \; \label{1}
 \end{eqnarray}
where $e$ is electric charge, $\Lambda$ is Cosmological constant, $M$ is Mass and $r$ is radial coordinates respectively.  


Throughout the paper, we will use the following expressions:
$Y_1=2Mr-r^2+\Lambda r^4-e^2,$ $Y_2=-Mr+\Lambda r^4+e^2,$ $Y_3=Mr+2\Lambda r^4-e^2,$ $Y_4=3Mr-4e^2,$ $Y_5=3Mr-2e^2,$ $Y_6=Mr-e^2,$ $Y_7=9Mr-10e^2$ and $Y_8=9Mr-8e^2.$

From the following calculations, we derive the non-vanishing Christoffel symbols of the second kind ($\Gamma^h_{ij}$), expressed as:
 \begin{eqnarray}\label{2}
 \begin{cases}

 	\Gamma^2_{11}=\frac{Y_1 Y_2}{r^5}, \ \ 
 	\Gamma^1_{12}=-\frac{2\Lambda r^3+M-r}{Y_1}, \\
 	\Gamma^3_{23}=\frac{1}{r}=\Gamma^4_{24}, \ \ 
 	\Gamma^2_{33}=\frac{Y_1}{r}=\frac{1}{\sin^2\theta}\Gamma^2_{44},\\
 	\Gamma^4_{34}=\cot\theta,\ \ \Gamma^3_{44}=-\cos\theta \sin\theta. \\
\end{cases}
  \end{eqnarray}

The following are the computed non-zero components of the Riemann curvature tensor and the Ricci tensor, which may potentially be non-vanishing
\begin{eqnarray}\label{3}
 \begin{cases}

 	R_{1212}=-\Lambda-\frac{2Mr-3e^2}{r^4}, \ \ 	
 	R_{1313}=\frac{Y_1 Y_2}{r^4}=\frac{1}{\sin^2\theta} R_{1414}, \\
 	R_{2323}=\frac{Y_2}{Y_1}=-\frac{1}{\sin^2\theta}R_{2424}, \  \
 	R_{3434}=(2Mr+\Lambda r^4-e^2) \sin^2\theta; \\
 \end{cases}
  \end{eqnarray}
 
 \begin{eqnarray}\label{4}
 \begin{cases}
 	Ric_{11}=-\frac{Y_1(3\Lambda-e^2)}{r^4}, \ \ 
 	Ric_{22}=\frac{3\Lambda-e^2}{r^2Y_1}, \\
 	Ric_{33}=\frac{3\Lambda-e^2}{r^2}=\frac{1}{\sin^2\theta}Ric_{44}. \ \ 
  \end{cases}
   \end{eqnarray}

  Finally, the scalar curvature is given by
  \begin{equation}
      \kappa=g_{ab}\,Ric^{ab}=-12\Lambda.\label{5}
  \end{equation}
  
%


Let $\mathcal{E}^{1} = R \cdot R,$ $\mathcal{E}^{2} =\nabla R,$ $\mathcal{E}^{3} = R \cdot Ric,$ $\mathcal{J}^1 = Q(g, R)$ and $\mathcal{J}^2 = Q(g, Ric).$ Taking into account their inherent symmetries, the non-vanishing components of $\mathcal{E}^{1},$ $\mathcal{E}^{2},$ $\mathcal{E}^{3},$ $\mathcal{J}^1,$ and $\mathcal{J}^2$ are listed below:

\begin{eqnarray}\label{RR}
\begin{cases}
	\mathcal{E}^1_{1223,13}=\frac{Y_2 Y_4}{r^6}=\mathcal{E}^1_{1224,14}=-\mathcal{E}^1_{1213,23}=-\frac{1}{\sin^2\theta}\mathcal{E}^1_{1214,24}, \\ 
	\mathcal{E}^1_{1434,13}=\frac{Y_1 Y_2 Y_5}{r^2}=-\mathcal{E}^1_{1334,14}, \\
	\mathcal{E}^1_{2434,23}=\frac{Y_2 Y_5\sin^2\theta}{r^2Y_1}=-\mathcal{E}^1_{2334,24}; 
	\end{cases}
	\end{eqnarray}

\begin{eqnarray}\label{R}
\begin{cases}
\mathcal{E}^{2}_{1212,2}=\frac{6(Mr-2e^2)}{r^5}=-\frac{r^4}{\sin^2\theta}\mathcal{E}^{2}_{2434,3}=-\frac{r^4}{2\sin^2\theta}\mathcal{E}^{2}_{2334,4},\\
\mathcal{E}^{2}_{1213,3}=\frac{Y_4 Y_1}{r^5}=\mathcal{E}^{2}_{1313,2}=\frac{1}{\sin^2\theta}\mathcal{E}^{2}_{1414,2}=\frac{1}{\sin^2\theta}\mathcal{E}^{2}_{1214,4},\\
\mathcal{E}^{2}_{2323,2}=-\frac{Y_4}{rY_1}=\frac{1}{\sin^2\theta}\mathcal{E}^{2}_{2424,2};\\
	\end{cases}
	\end{eqnarray}

\begin{eqnarray}\label{RS}
\begin{cases}
\mathcal{E}^{3}_{1313}=\frac{2e^2 Y_1 Y_2}{r^8}=\frac{1}{\sin^2\theta}\mathcal{E}^{3}_{1414},\\
\mathcal{E}^{3}_{2323}=-\frac{2e^2 Y_2}{r^4 Y_1}=\frac{1}{\sin^2\theta}\mathcal{E}^{3}_{2424};\\
\end{cases}
	\end{eqnarray}
\begin{eqnarray}\label{Q(g,R)}
\begin{cases}	
\mathcal{J}^1_{1223,13}=-\frac{Y_4}{r}=\frac{1}{\sin^2\theta}\mathcal{J}^1_{1224,14}=-\mathcal{J}^1_{1213,23}=-\frac{1}{\sin^2\theta}\mathcal{J}^1_{1214,24},\\
	\mathcal{J}^1_{1434,13}=\frac{Y_1 Y_5 \sin^2\theta}{r^2}=-\mathcal{J}^1_{1334,14}, \\
	\mathcal{J}^1_{2434,23}=-\frac{r^2 Y_5\sin^2\theta}{Y_1}=-\mathcal{J}^1_{2334,24}; \\
	\end{cases}
		\end{eqnarray}

\begin{eqnarray}\label{Q(g,S)}
\begin{cases}	
	\mathcal{J}^2_{1313}=-\frac{2e^2 Y_1}{r^4}=\frac{1}{\sin^2\theta}\mathcal{J}^2_{1414}, \\ 
	\mathcal{J}^2_{2323}=\frac{2e^2}{Y_1}=\frac{1}{\sin^2\theta}\mathcal{J}^2_{2424}.
	\end{cases}
		\end{eqnarray}

From the expressions in equations \eqref{RR}, \eqref{R},  \eqref{Q(g,R)},\eqref{RS} and \eqref{Q(g,S)}, We access the following results:

 \begin{pr}\label{pr1}
	The RNdS spacetime satisfies the following curvature conditions:
	\begin{enumerate}[label=(\roman*)]
	 \item $\ R\cdot R=-\frac{e^2-Mr+\Lambda r^4}{r^4}Q(g,R),$ i.e., the nature of the spacetime is pseudosymmetric\cite{Kowa06},

\item  $\ R\cdot Ric=-\frac{e^2-Mr+\Lambda r^4}{r^4}Q(g, Ric),$ i.e., the nature of the spacetime is Ricci peudosymmetric,

\item the general form of $R$-compatible tensor is given by	
	$$
	\left(
	\begin{array}{cccc}
		\mathscr{Z}_{11} &\mathscr{Z}_{12} & 0 & 0 \\
		\mathscr{Z}_{12} & \mathscr{Z}_{22} & 0 & 0 \\
		0 & 0 & \mathscr{Z}_{33} & \mathscr{Z}_{34} \\
		0 & 0 & \mathscr{Z}_{34} & \mathscr{Z}_{44}
	\end{array}
	\right),
	$$
	where $\mathscr{Z}_{ij}$ are arbitrary scalars.  
	 
	 \end{enumerate}
	
\end{pr}
Taking into account the symmetry properties of the conformal curvature tensor $C_{hijk}$, its non-vanishing components are given below.

$$\begin{array}{c}
	C_{1212}=-\frac{6 Y_6}{r^4}=\frac{1}{\sin^2\theta}C_{3434},\\ 
	C_{1313}=-\frac{Y_1 Y_6}{r^4}=-\frac{1}{\sin^2\theta}C_{1414}, \\
	C_{2323}=\frac{Y_6}{Y_1}=\frac{1}{\sin^2\theta}C_{2424}. \ \ 
\end{array}$$

The non-zero components of covariant derivatives of conformal curvature tensor are:
$$\begin{array}{c}
C_{1212,2}=\frac{2Y_4}{r^5},\\
C_{1213,3}=\frac{3Y_1 Y_5}{r^5}=\frac{1}{\sin^2\theta}C_{1214,4},\\
C_{1313,2}=\frac{Y_1 Y_4}{r^5}=C_{1414,2},\\
C_{2323,2}=-\frac{Y_4}{rY_1}=\frac{1}{\sin^2\theta}C_{2424,2},\\
C_{2334,4}=\frac{3 Y_5\sin^2\theta}{r}=-C_{2434,3},\\
C_{3434,2}=-\frac{2Y_4\sin^2\theta}{r};
\end{array}$$

Let  $\mathcal{E}^{4}=R\cdot C$, $\mathcal{E}^{5}=C\cdot R$, $\mathcal{E}^{6}=C\cdot C$ , $\mathcal{E}^{7}=C\cdot Ric,$ $\mathcal{J}^3=Q(g,C),$ $\mathcal{J}^4=Q(Ric,R)$ and $\mathcal{J}^5=Q(Ric,C)$.  Considering their symmetry properties, the non-zero components of $\mathcal{E}^{4}$, $\mathcal{E}^{5}$, $\mathcal{E}^{6}$, $\mathcal{E}^{7}$, $\mathcal{J}^3$, $\mathcal{J}^4$ and $\mathcal{J}^5$ are presented below:

\begin{eqnarray}\label{RC}
\begin{cases}
		\mathcal{E}^{4}_{1223,13}=-\frac{3Y_2 Y_6}{r^6}=\frac{1}{\sin^2\theta}\mathcal{E}^{4}_{1224,14}=-\mathcal{E}^{4}_{1213,23}=-\frac{1}{\sin^2\theta}\mathcal{E}^{4}_{1214,24}, \\ 
		\mathcal{E}^{4}_{1434,13}=-\frac{Y_1 Y_2 Y_6}{r^6}=-\mathcal{E}^{4}_{1334,14},\\
		\mathcal{E}^{4}_{2434,23}=-\frac{3Y_2 Y_6\sin^ 2\theta}{r^2Y_1}=-\mathcal{E}^{4}_{2334,24}; \\
	
	\end{cases}
			\end{eqnarray}

\begin{eqnarray}\label{CR}
\begin{cases}
	\mathcal{E}^{5}_{1223,13}=-\frac{Y_4Y_5}{r^6}=\frac{1}{\sin^2\theta}\mathcal{E}^{5}_{1224,14}=-\mathcal{E}^{5}_{1213,23}=-\frac{1}{\sin^2\theta}\mathcal{E}^{5}_{1214,24}, \\
	\mathcal{E}^{5}_{1434,13}=\frac{Y_1 Y_5 Y_6\sin^2\theta}{r^6}=-\mathcal{E}^{5}_{1334,14}, \\
	\mathcal{E}^{5}_{2434,23}=-\frac{Y_5Y_6\sin^2\theta}{r^2Y_1}=-\mathcal{E}^{5}_{2334,24}; 
	\end{cases}
			\end{eqnarray}

\begin{eqnarray}\label{CC}
\begin{cases}
\mathcal{E}^6_{1223,13}=-\frac{3Y_6^2}{r^6}=\frac{1}{\sin^2\theta}\mathcal{E}^{6}_{1224,14}=-\mathcal{E}^{6}_{1213,23}=\frac{1}{\sin^2\theta}\mathcal{E}^{6}_{1214,24}, \\ 
	\mathcal{E}^{6}_{1434,13}=\frac{3Y_1 Y_6\sin^2\theta}{r^6}=-\mathcal{E}^{6}_{1334,14}, \\
	\mathcal{E}^{6}_{2434,23}=\frac{3Y_6\sin^2\theta}{r^2Y_1}=-\mathcal{E}^{6}_{2334,24};
\end{cases}
		\end{eqnarray}
		\begin{eqnarray}\label{CS}
		\begin{cases}
		\mathcal{E}^{7}_{13,13}=-\frac{2e^2 Y_1 Y_6}{r^8}=\frac{1}{\sin^2\theta}\mathcal{E}^{7}_{14,14},\\
		\mathcal{E}^{7}_{23,23}=-\frac{2e^2 Y_6}{r^4 Y_1}=\frac{1}{\sin^2\theta}\mathcal{E}^{7}_{24,24};\\
		\end{cases}
			\end{eqnarray}
\begin{eqnarray}\label{Q(g,C)}
\begin{cases}
\mathcal{J}^3_{1223,13}=-\frac{3Y_6}{r^2}=\frac{1}{\sin^2\theta}\mathcal{J}^3_{1224,14}=-\mathcal{J}^3_{1213,23}=-\frac{1}{\sin^2\theta}\mathcal{J}^3_{1214,24},\\	
	\mathcal{J}^3_{1434,13}=\frac{3Y_1 Y_6 \sin^2\theta}{r^2}=-\mathcal{J}^3_{1334,14}, \\
	\mathcal{J}^3_{2434,23}=-\frac{3r^2Y_6\sin^2\theta}{Y_1}=-\mathcal{J}^3_{2334,24}; \\
\end{cases}
		\end{eqnarray}

\begin{eqnarray}\label{Q(Ric,R)}
\begin{cases}
\mathcal{J}^4_{1223,13}=-\frac{9M\Lambda r^5+r(M-10\Lambda r^3)e^2-2e^4}{r^6}=\frac{1}{\sin^2\theta}\mathcal{J}^4_{1224,14}=-\mathcal{J}^4_{1213,23}=\frac{1}{\sin^2\theta}\mathcal{J}^4_{1214,24}, \\ 
	\mathcal{J}^4_{1434,13}=\frac{Y_1 (9M\Lambda r^4+r(M-8\Lambda r^3)e^2)\sin^2\theta}{r^5}=-\mathcal{J}^4_{1334,14}, \\
		\mathcal{J}^4_{2434,23}=-\frac{9M\Lambda r^4+r(M-8\Lambda r^3)e^2\sin^2\theta}{rY_1}=-\mathcal{J}^4_{2334,24}; 	
\end{cases}
		\end{eqnarray}
\begin{eqnarray}\label{Q(Ric,C)}
\begin{cases}
\mathcal{J}^5_{1223,13}=\frac{(9\Lambda r^4+e^2)Y_6}{r^6}=\frac{1}{\sin^2\theta}\mathcal{J}^5_{1224,14}=-\mathcal{J}^5_{1213,23}=\frac{1}{\sin^2\theta}\mathcal{J}^5_{1214,24}, \\ 
	\mathcal{J}^5_{1434,13}=\frac{(9\Lambda r^4-e^2)Y_1Y_6\sin^2\theta}{r^6}=-\mathcal{J}^5_{1334,14}, \\
		\mathcal{J}^5_{2434,23}=\frac{(9\Lambda r^4-e^2)Y_6\sin^2\theta}{r^2Y_1}=-\mathcal{J}^5_{2334,24}; 	
\end{cases}
		\end{eqnarray}		
		
		Using the expressions given in equations \eqref{RC}, \eqref{CR}, \eqref{CC}, \eqref{Q(g,C)}, \eqref{Q(g,C)}, \eqref{Q(Ric,R)} and \eqref{Q(Ric,C)}, we obtain the following outcomes:


 \begin{pr}\label{pr2} The RNdS spacetime satisfies the specified curvature conditions: 
 \begin{enumerate}[label=(\roman*)]
 				
 \item	$\ R\cdot C=-\frac{e^2-Mr+\Lambda r^4}{r^4}Q(g,C)$
 	 i.e., the spacetime is conformally pseudosymmetric.
 
\item $C\cdot Ric=\frac{Mr-e^2}{r^4}Q(g,Ric)$	i.e., Ricci pseudosymmetric due to conformal tensor.
\item	$\ C\cdot C=\frac{Mr-e^2}{r^4}Q(g,C)\ \mbox{and} \ C\cdot R=\frac{Mr-e^2}{r^4}Q(g,R)$
	 i.e., the spacetime is pseudosymmetric due to conformal tensor and also Ricci generalized conformal peudosymmetric manifold.
	 
	  \item $\ R\cdot R-Q(Ric,R)=\frac{2e^4-6e^2Mr+3M^2r^2-6\Lambda e^2r^4+6M\Lambda r^5)}{3r^4(Mr-e^2)}Q(g,C).$
\item $C \cdot R - R \cdot C= \mathcal{L}_1Q(g, R) +\mathcal{L}_2 Q(Ric,R),$\\
where $\mathcal{L}_1=\frac{Y_6}{r^4}\left \lbrace 1-\frac{3Y_2^2}{3Mr(Mr-2e^2)-6\Lambda r^4Y_6+2e^4} \right\rbrace$ and
$\mathcal{L}_2=\left \lbrace\frac{3Y_6Y_2}{3Mr(Mr-2e^2)-6\Lambda r^4Y_6+2e^4}\right\rbrace.$
\item $C \cdot R - R \cdot C=4\Lambda Q(g, C) +Q(Ric,C).$	  
	  \end{enumerate}	
\end{pr}
Taking into account the symmetry characteristics of the projective curvature tensor $P_{hijk}$, its non-vanishing components are listed below.

$$\begin{array}{c}
	P_{1212}=-\frac{2Y_4}{3r^4}=-P_{1221},  \\
	P_{1313}=-\frac{Y_1 Y_4}{3r^4}=\frac{1}{\sin^2\theta}P_{1414}, \\
    P_{1331}=\frac{Y_1 Y_5}{3r^4}=\frac{1}{\sin^2\theta}P_{1441},\\
	P_{2323}=\frac{Y_4}{3Y_1}=\frac{1}{\sin^2\theta}P_{2424},\\
	P_{2332}=-\frac{Y_5}{3Y_1}=-\frac{1}{\sin^2\theta}P_{2442}, \\
	P_{3434}=\frac{2Y_5}{3}=-P_{3443}.
\end{array}$$

Let $\mathcal{E}^{8}=R\cdot P,$ $\mathcal{E}^{9}=C\cdot P,$ $\mathcal{E}^{10}=P\cdot Ric$ and  $\mathcal{J}^6=Q(g,p)$.  Taking their symmetries into account, the non-zero components of $\mathcal{E}^{8}$, $\mathcal{E}^{9}$, $\mathcal{E}^{10}$ and $\mathcal{J}^6$ are given as follows.


\begin{eqnarray}\label{RP}
\begin{cases}
    \mathcal{E}^{8}_{1223,13}=-\frac{Y_2 Y_4}{r^6}=-\mathcal{E}^{8}_{2321,13}=\frac{1}{\sin^2\theta}\mathcal{E}^{8}_{1224,14}=-\frac{1}{\sin^2\theta}\mathcal{E}^{8}_{2421,14}\\=-\mathcal{E}^{8}_{1213,23}=-\mathcal{E}^{8}_{1312,23}=-\frac{1}{\sin^2\theta}\mathcal{E}^{8}_{1214,24}=-\frac{1}{\sin^2\theta}\mathcal{E}^{8}_{1412,24},\\\mathcal{E}^{8}_{1232,13}-\frac{Y_2 Y_7}{r^6}=-\mathcal{E}^{8}_{2312,13}=-\frac{1}{\sin^2\theta}\mathcal{E}^{8}_{1242,14}=\frac{1}{\sin^2\theta}\mathcal{E}^{8}_{2412,14}\\=\mathcal{E}^{8}_{1231,23}=\mathcal{E}^{8}_{1321,23}=\frac{1}{\sin^2\theta}\mathcal{E}^{8}_{1241,24}=\frac{1}{\sin^2\theta}\mathcal{E}^{8}_{1421,24},\\
    \mathcal{E}^{8}_{1311,13}=\frac{2e^2 Y_1 Y_2}{3r^{10}}=\frac{1}{\sin^2\theta}\mathcal{E}^{8}_{1411,14},\\
    \mathcal{E}^{8}_{1333,13}=\frac{2e^2 Y_1 Y_2}{3r^6}=\frac{1}{\sin^2\theta}\mathcal{E}^{8}_{2444,14},\\
    \mathcal{E}^{8}_{1434,13}=\frac{Y_1 Y_2 Y_8\sin^2\theta}{3r^6} =\mathcal{E}^{8}_{3414,13}=-\mathcal{E}^{8}_{1343,14}=-\mathcal{E}^{8}_{3413,14},\\
    
 \mathcal{E}^{8}_{1443,13}=\frac{Y_1Y_2Y_5\sin^2\theta}{r^6}=-\mathcal{E}^{8}_{3441,13}=-\mathcal{E}^{8}_{1334,14}=\mathcal{E}^{8}_{3431,14},\\
 \mathcal{E}^{8}_{2322,23}=-\frac{2e^2Y_2}{3r^2Y_1^2}=\frac{1}{\sin^2\theta}\mathcal{E}^{8}_{2422,24},\\
     \mathcal{E}^{8}_{2333,23}=-\frac{2e^2Y_2}{3r^2Y_1} =\frac{1}{\sin^2\theta}\mathcal{E}^{8}_{2444,24},\\
    \mathcal{E}^{8}_{2434,23}=-\frac{Y_2Y_8\sin^2\theta}{3r^2Y_1}=\mathcal{E}^{8}_{3424,23}=\mathcal{E}^{8}_{2343,24}=\mathcal{E}^{8}_{3423,24}\\
    \mathcal{E}^{8}_{2443,23}=-\frac{Y_2Y_5\sin^2\theta}{r^2Y_1}=-\mathcal{E}^{8}_{3442,23} =\mathcal{E}^{8}_{2334,24}=\mathcal{E}^{8}_{3432,24};\\
\end{cases}
		\end{eqnarray}

\begin{eqnarray}\label{CP}
\begin{cases}
 \mathcal{E}^{9}_{1223,13}=-\frac{Y_4 Y_6}{r^6}=-\mathcal{E}^{9}_{2321,13}=\frac{1}{\sin^2\theta}\mathcal{E}^{9}_{1224,14}=-\frac{1}{\sin^2\theta}\mathcal{E}^{9}_{2421,14}\\=-\mathcal{E}^{9}_{1213,23}=-\mathcal{E}^{9}_{1312,23}=-\frac{1}{\sin^2\theta}\mathcal{E}^{9}_{1214,24}=-\frac{1}{\sin^2\theta}\mathcal{E}^{9}_{1412,24},\\
 \mathcal{E}^{9}_{1232,13}-\frac{Y_6 Y_7}{3r^6}=-\mathcal{E}^{9}_{2312,13}=-\frac{1}{\sin^2\theta}\mathcal{E}^{9}_{1242,14}=\frac{1}{\sin^2\theta}\mathcal{E}^{9}_{2412,14}\\=\mathcal{E}^{9}_{1231,23}=\mathcal{E}^{9}_{1321,23}=\frac{1}{\sin^2\theta}\mathcal{E}^{9}_{1241,24}=\frac{1}{\sin^2\theta}\mathcal{E}^{9}_{1421,24},\\
    \mathcal{E}^{9}_{1311,13}=-\frac{2e^2 Y_1 Y_6}{3r^{10}}=\frac{1}{\sin^2\theta}\mathcal{E}^{9}_{1411,14},\\
    \mathcal{E}^{9}_{1333,13}=-\frac{2e^2 Y_1 Y_6}{3r^6}=\frac{1}{\sin^2\theta}\mathcal{E}^{9}_{2444,14},\\
     \mathcal{E}^{9}_{1434,13}=\frac{Y_1 Y_6 Y_8\sin^2\theta}{3r^6} =\mathcal{E}^{9}_{3414,13}=-\mathcal{E}^{9}_{1343,14}=-\mathcal{E}^{9}_{3413,14},\\
     
  \mathcal{E}^{9}_{1443,13}=\frac{Y_1Y_2Y_5\sin^2\theta}{r^6}=-\mathcal{E}^{9}_{3441,13}=-\mathcal{E}^{9}_{1334,14}=\mathcal{E}^{9}_{3431,14},\\
  \mathcal{E}^{9}_{2322,23}=-\frac{2e^2Y_6}{3r^2Y_1^2}=\frac{1}{\sin^2\theta}\mathcal{E}^{9}_{2422,24},\\
      \mathcal{E}^{9}_{2333,23}=\frac{2e^2Y_6}{3r^2Y_1} =\frac{1}{\sin^2\theta}\mathcal{E}^{9}_{2444,24},\\
     \mathcal{E}^{9}_{2434,23}=-\frac{Y_6Y_8\sin^2\theta}{3r^2Y_1}=\mathcal{E}^{9}_{3424,23}=\mathcal{E}^{9}_{2343,24}=\mathcal{E}^{9}_{3423,24}\\
     \mathcal{E}^{9}_{2443,23}=-\frac{Y_5Y_6\sin^2\theta}{r^2Y_1}=-\mathcal{E}^{9}_{3442,23} =\mathcal{E}^{9}_{2334,24}=\mathcal{E}^{9}_{3432,24};\\
\end{cases}
		\end{eqnarray}


\begin{eqnarray}\label{PS}
\begin{cases}
 \mathcal{E}^{10}_{13,13}=\frac{2e^2Y_1Y_2}{r^8}=\frac{1}{\sin^2\theta}\mathcal{E}^{10}_{14,14},\\
 \mathcal{E}^{10}_{23,23}=-\frac{2e^2Y_2}{r^4Y_1}=\frac{1}{\sin^2\theta}\mathcal{E}^{10}_{24,24},\\
 \mathcal{E}^{10}_{13,31}=-\frac{2e^2Y_1Y_2}{r^8}=-\frac{1}{\sin^2\theta}\mathcal{E}^{10}_{14,41},\\
      \mathcal{E}^{10}_{23,32}=\frac{2e^2Y_2}{r^4Y_1}=\frac{1}{\sin^2\theta}\mathcal{E}^{10}_{24,42};\\
\end{cases}
		\end{eqnarray}

\begin{eqnarray}\label{Q(g,P)}
\begin{cases}
  \mathcal{J}^6_{1223,13}=-\frac{Y_4}{r^2}=-\mathcal{J}^6_{2321,13}=\frac{1}{\sin^2\theta}\mathcal{J}^6_{1224,14}=-\frac{1}{\sin^2\theta}\mathcal{J}^6_{2421,14}\\=-\mathcal{J}^6_{1213,23}=-\mathcal{J}^6_{1312,23}=-\frac{1}{\sin^2\theta}\mathcal{J}^6_{1214,24}=-\frac{1}{\sin^2\theta}\mathcal{J}^6_{1412,24}
  \\\mathcal{J}^6_{1232,13}-\frac{Y_7}{3r^2}=-\mathcal{J}^6_{2312,13}=-\frac{1}{\sin^2\theta}\mathcal{J}^6_{1242,14}=\frac{1}{\sin^2\theta}\mathcal{J}^6_{2412,14}\\=\mathcal{J}^6_{1231,23}=\mathcal{J}^6_{1321,23}=\frac{1}{\sin^2\theta}\mathcal{J}^6_{1241,24}=\frac{1}{\sin^2\theta}\mathcal{J}^6_{1421,24},\\
     \mathcal{J}^6_{1311,13}=\frac{2e^2 Y_1^2}{3r^6}=\frac{1}{\sin^2\theta}\mathcal{J}^6_{1411,14},\\
     \mathcal{J}^6_{1333,13}=\frac{2e^2 Y_1}{3r^2}=\frac{1}{\sin^2\theta}\mathcal{J}^6_{2444,14},\\
     \mathcal{J}^6_{1434,13}=\frac{Y_1 Y_8\sin^2\theta}{3r^2} =\mathcal{J}^6_{3414,13}=-\mathcal{J}^6_{1343,14}=-\mathcal{J}^6_{3413,14},\\
     
  \mathcal{J}^6_{1443,13}=\frac{Y_1Y_5\sin^2\theta}{r^2}=-\mathcal{J}^6_{3441,13}=-\mathcal{J}^6_{1334,14}=\mathcal{J}^6_{3431,14},\\
  \mathcal{J}^6_{2322,23}=-\frac{2r^2e^2}{32Y_1^2}=\frac{1}{\sin^2\theta}\mathcal{J}^6_{2422,24},\\
      \mathcal{J}^6_{2333,23}=-\frac{2r^2e^2}{3Y_1} =\frac{1}{\sin^2\theta}\mathcal{J}^6_{2444,24},\\
     \mathcal{J}^6_{2434,23}=-\frac{r^2Y_8\sin^2\theta}{3Y_1}=\mathcal{J}^6_{3424,23}=\mathcal{J}^6_{2343,24}=\mathcal{J}^6_{3423,24}\\
     \mathcal{J}^6_{2443,23}=\frac{r^2Y_5\sin^2\theta}{Y_1}=-\mathcal{J}^6_{3442,23} =\mathcal{J}^6_{2334,24}=\mathcal{J}^6_{3432,24};\\
\end{cases}
		\end{eqnarray}


		
	From the expressions given in equations \eqref{RP}, \eqref{CP}, \eqref{PS}, and \eqref{Q(g,P)}, the following results are obtained.
 \begin{pr}\label{pr3}
The RNdS spacetime satisfies the given curvature relations. 
\begin{enumerate}[label=(\roman*)]
	\item  $R\cdot P=-\frac{e^2-Mr+\Lambda r^4}{r^4}Q(g,P),$

	\item  $C\cdot P=-\frac{e^2-Mr}{r^4}Q(g,P),$

 
	\item $P\cdot Ric=-\frac{e^2-Mr+\Lambda r^4}{r^4}Q(g,Ric).$
	
		\end{enumerate}
\end{pr}

Taking into account the symmetry properties of the concircular curvature tensor $W_{hijk}$, its non-zero components are listed below:

$$\begin{array}{c}
W_{1212}=\frac{-2Mr+3e^2}{r^4},
	W_{1313}=-\frac{Y_1 Y_6}{r^4}=\frac{1}{\sin^2\theta}W_{1414}, \\
	W_{2323}=\frac{Y_6}{Y_1}=\frac{1}{\sin^2\theta}W_{2424}, 
	W_{3434}=(2Mr-e^2)\sin^2\theta. \\ 
\end{array}$$

Let $\mathcal{I}^1=R \cdot W,$ $\mathcal{I}^2=C \cdot W,$ $\mathcal{I}^3=W \cdot R,$ $\mathcal{I}^{4}=W \cdot Ric,$ $\mathcal{I}^{5}=W \cdot C,$ $\mathcal{I}^{6}=W \cdot P,$ $\mathcal{I}^{7}=W \cdot W$ and $\mathcal{J}^7=Q(g,W).$  Considering their symmetry properties, the non-vanishing components of $\mathcal{I}^1$, $\mathcal{I}^2$, $\mathcal{I}^3$, $\mathcal{I}^4$, $\mathcal{I}^5$, $\mathcal{I}^6,$ $\mathcal{I}^{7}$ and $\mathcal{J}^7$ are presented below.

 

\begin{eqnarray}\label{RW}
\begin{cases}
	\mathcal{I}^1_{1223,13}=\frac{Y_2Y_4}{r^6}=\frac{1}{\sin^2\theta}\mathcal{I}^1_{1224,14}= -\mathcal{I}^1_{1213,23} =-\frac{1}{\sin^2\theta}\mathcal{I}^1_{1214,24}, \\
	\mathcal{I}^1_{1434,13}=-\frac{Y_1Y_2Y_5\sin^2\theta}{r^6}=-\mathcal{I}^1_{1334,14}, \\
	\mathcal{I}^1_{2434,23}=\frac{Y_2Y_5\sin^2\theta}{r^2Y_1}=-\mathcal{I}^1_{2334,24};
\end{cases}
		\end{eqnarray}

\begin{eqnarray}\label{CW}
\begin{cases} \mathcal {I}^2_{1223,13}=-\frac{Y_4Y_6}{r^6}=\frac{1}{\sin^2\theta}\mathcal{I}^2_{1224,14}=-\mathcal {I}^2_{1213,23}=\frac{1}{\sin^2\theta}\mathcal {I}^2_{1214,24}, \\ 
	\mathcal {I}^2_{1434,13}=\frac{(Mr-r^2+\Lambda R^4-e^2)Y_5Y_6\sin^2\theta}{r^6}=-\mathcal {I}^2_{1334,14}, \\
	\mathcal {I}^2_{2434,23}=-\frac{Y_5Y_6}{r^2Y_1}=-\mathcal {I}^2_{2334,24};
\end{cases}
		\end{eqnarray}


\begin{eqnarray}\label{WR}
\begin{cases}
\mathcal {I}^3_{1223,13}=-\frac{Y_4Y_6}{r^6}=\frac{1}{\sin^2\theta}\mathcal{I}^3_{1224,14}=-\mathcal {I}^3_{1213,23}=\frac{1}{\sin^2\theta}\mathcal{I}^3_{1214,24}, \\ 
	\mathcal{I}^3_{1434,13}=\frac{(Mr-r^2+\Lambda R^4-e^2)Y_5Y_6\sin^2\theta}{r^6}=-\mathcal{I}^3_{1334,14}, \\
	\mathcal{I}^3_{2434,23}=-\frac{Y_5Y_6}{r^2Y_1}=-\mathcal{I}^3_{2334,24};
\end{cases}
		\end{eqnarray}
		\begin{eqnarray}\label{WS}
		\begin{cases}
		\mathcal{I}^{4}_{13,13}=-\frac{2e^2 Y_1 Y_6}{r^8}=\frac{1}{\sin^2\theta}\mathcal{I}^{4}_{14,14},\\
		\mathcal{I}^{4}_{23,23}=-\frac{2e^2 Y_6}{r^4 Y_1}=\frac{1}{\sin^2\theta}\mathcal{I}^{4}_{24,24};\\
		\end{cases}
			\end{eqnarray}		
		



\begin{eqnarray}\label{WC}
\begin{cases}
\mathcal{I}^{5}_{1223,13}=-\frac{3Y_6^2}{r^6}=\frac{1}{\sin^2\theta}\mathcal{I}^{5}_{1224,14}=-\mathcal{I}^{5}_{1213,23}=\frac{1}{\sin^2\theta}\mathcal{I}^{5}_{1214,24}, \\ 
	\mathcal{I}^{5}_{1434,13}=\frac{3Y_1Y_6^2\sin^2\theta}{r^6}=-\mathcal{I}^{5}_{1334,14}, \\
	\mathcal{I}^{5}_{2434,23}=-\frac{3Y_6^2\sin^2\theta}{r^2Y_1}=-\mathcal{I}^{5}_{2334,24};
\end{cases}
		\end{eqnarray}
%

%

\begin{eqnarray}\label{WP}
\begin{cases}
    \mathcal{I}^{6}_{1223,13}=-\frac{Y_4Y_6}{r^6}=\frac{1}{\sin^2\theta}\mathcal{I}^{6}_{1224,14}=-\mathcal{I}^{6}_{1213,23}=-\frac{1}{\sin^2\theta}\mathcal{I}^{6}_{1214,24},
    \\\mathcal{I}^{6}_{1232,13}=\frac{Y_6 Y_7}{3r^6}=-\frac{1}{\sin^2\theta}\mathcal{I}^{6}_{1242,14}=\mathcal{I}^{6}_{1231,23}=\mathcal{I}^{6}_{1321,23}\\=-\frac{1}{\sin^2\theta}\mathcal{I}^{6}_{1241,24}=-\frac{1}{\sin^2\theta}\mathcal{I}^{6}_{1421,24},\\
    \mathcal{I}^{6}_{1333,13}=-\frac{2e^2 Y_1 Y_6}{3r^6}=\frac{1}{\sin^2\theta}\mathcal{I}^{6}_{1444,14},\\
    \mathcal{I}^{6}_{1434,13}=\frac{Y_1 Y_6 Y_8\sin^2\theta}{3r^6} =-\mathcal{I}^{6}_{1343,14},\\
    
 \mathcal{I}^{6}_{1443,13}=\frac{Y_1Y_5Y_6\sin^2\theta}{r^6}=-\mathcal{I}^{6}_{3441,13}=-\mathcal{I}^{6}_{1334,14},\\
     \mathcal{I}^{6}_{2333,23}=\frac{2e^2Y_6}{3r^2Y_1} =\frac{1}{\sin^2\theta}\mathcal{I}^{6}_{2444,24},\\
    \mathcal{I}^{6}_{2434,23}=-\frac{Y_6Y_8\sin^2\theta}{3r^2Y_1}=\mathcal{I}^{6}_{2343,24}\\
    \mathcal{I}^{6}_{2443,23}=\frac{Y_5Y_6\sin^2\theta}{r^2Y_1}=-\mathcal{I}^{6}_{3442,23} =\mathcal{I}^{6}_{2334,24};\\
\end{cases}
		\end{eqnarray}


\begin{eqnarray}\label{WW}
\begin{cases}
\mathcal{I}^{7}_{1223,13}=-\frac{Y_4Y_6}{r^6}=\frac{1}{\sin^2\theta}\mathcal{I}^{7}_{1224,14}=-\mathcal{I}^{7}_{1213,23}=\frac{1}{\sin^2\theta}\mathcal{I}^{7}_{1214,24}, \\ 
	\mathcal{I}^{7}_{1434,13}=\frac{(Mr-r^2+\Lambda r^4-e^2)Y_5Y_6\sin^2\theta}{r^6}=-\mathcal{I}^{7}_{1334,14}, \\
	\mathcal{I}^{7}_{2434,23}=-\frac{Y_5Y_6}{r^2Y_1}=-\mathcal{I}^{7}_{2334,24};
\end{cases}
		\end{eqnarray}

\begin{eqnarray}\label{Q(g,W)}
\begin{cases}
	\mathcal{J}^7_{1223,13}=-\frac{Y_4}{r^2}=\frac{1}{\sin^2\theta}\mathcal{J}^7_{1224,14}= -\mathcal{J}^7_{1213,23} =-\frac{1}{\sin^2\theta}\mathcal{J}^7_{1214,24}, \\
		\mathcal{J}^7_{1434,13}=\frac{Y_1Y_5\sin^2\theta}{r^2}=-\mathcal{J}^7_{1334,14}, \\
		\mathcal{J}^7_{2434,23}=-\frac{r^2Y_4\sin^2\theta}{Y_1}=-\mathcal{J}^7_{2334,24}. 
\end{cases}
		\end{eqnarray}

Referring to equations \eqref{RW}, \eqref{CW}, \eqref{WR}, \eqref{WS}, \eqref{WC}, \eqref{WP}, and \eqref{Q(g,W)}, the following results are obtained.

 \begin{pr}\label{pr4}
	The RNdS spacetime satisfies the following curvature conditions of pseudosymmetric type:
		\begin{enumerate}[label=(\roman*)]
				
	 	\item  $R\cdot W=-\frac{e^2-Mr+\Lambda r^4}{r^4}Q(g,W),$
	
		\item  $C\cdot W=-\frac{e^2-Mr}{r^4}Q(g,W),$
	
		
			\item  $W\cdot R=-\frac{e^2-Mr}{r^4}Q(g,R)\ \mbox{and} \ W\cdot Ric=-\frac{e^2-Mr}{r^4}Q(g,Ric),$	
		\item  $W\cdot C=-\frac{e^2-Mr}{r^4}Q(g,C)\ \mbox{and} \ W\cdot P=-\frac{e^2-Mr}{r^4}Q(g,P),$
		
			\item  $W\cdot W=-\frac{e^2-Mr}{r^4}Q(g,W).$		
	 
	 \end{enumerate}	
\end{pr}
Taking into account its symmetry properties, the conharmonic curvature tensor $K_{hijk}$ has the following non-zero components.
	
$$\begin{array}{c}
	K_{1212}=\frac{2Y_2}{r^4},\ \ 
	K_{1313}=\frac{Y_1Y_2}{r^4}=\frac{1}{\sin^2\theta}K_{1414}, \\
K_{2323}=-\frac{Y_2}{Y_1}=\frac{1}{\sin^2\theta}K_{2424}, \ \ 
	K_{3434}=-Y_2\sin^2\theta.
\end{array}$$

Let $\mathcal{Y}^{1}=R \cdot K,$  $\mathcal{Y}^{2}=C \cdot K,$ $\mathcal{Y}^{3}=K \cdot Ric,$ $\mathcal{Y}^{4}=W \cdot K,$ $\mathcal{Y}^{5}=K \cdot R,$ $\mathcal{Y}^{6}=K \cdot C,$ $\mathcal{Y}^{7}=K \cdot P,$ $\mathcal{Y}^{8}=K \cdot W,$ $\mathcal{Y}^{9}= \ K \cdot K$ and $\mathcal{J}^8=Q(g,K).$ Taking their symmetries into consideration, the non-zero components of $\mathcal{Y}^{1}$, $\mathcal{Y}^{2}$, $\mathcal{Y}^{3}$, $\mathcal{Y}^{4}$, $\mathcal{Y}^{5},$ $\mathcal{Y}^{6},$ $\mathcal{Y}^{7},$ $\mathcal{Y}^{8},$ $\mathcal{Y}^{9}$ and $\mathcal{J}^8$ are determined as follows.

\begin{eqnarray}\label{RK}
\begin{cases}
	\mathcal{Y}^{1}_{1223,13}=-\frac{Y_2Y_6}{r^6}=\frac{1}{\sin^2\theta}\mathcal{Y}^{1}_{1224,14}=-\mathcal{Y}^{1}_{1213,23}=-\frac{1}{\sin^2\theta}\mathcal{Y}^{1}_{1214,24}, \\ 
	\mathcal{Y}^{1}_{1434,13}=\frac{3Y_1Y_2Y_6}{r^6}=-\mathcal{Y}^{1}_{1334,14},\\
	\mathcal{Y}^{1}_{2434,23}=\frac{3Y_2Y_6}{r^2Y_1}=-\mathcal{Y}^{1}_{2334,24}; 
\end{cases}
		\end{eqnarray}

\begin{eqnarray}\label{CK}
\begin{cases}
	\mathcal{Y}^{2}_{1223,13}=-\frac{3Y_6^2}{r^6}=\frac{1}{\sin^2\theta}\mathcal{Y}^{2}_{1224,14}=-\mathcal{Y}^{2}_{1213,23}=-\frac{1}{\sin^2\theta}\mathcal{Y}^{2}_{1214,24}, \\
	\mathcal{Y}^{2}_{1434,13}=\frac{3Y_1Y_6^2\sin^2\theta}{r^6}=-\mathcal{Y}^{2}_{1334,14}, \\
	\mathcal{Y}^{2}_{2434,23}=-\frac{3Y_6^2\sin^2\theta}{r^2Y_1}=-\mathcal{Y}^{2}_{2334,24}; 
\end{cases}
		\end{eqnarray}
	
	

\begin{eqnarray}\label{KS}
\begin{cases}
 \mathcal{Y}^{3}_{13,13}=\frac{2e^2Y_1Y_2}{r^8}=\frac{1}{\sin^2\theta}\mathcal{Y}^{3}_{14,14},\\
 \mathcal{Y}^{3}_{23,23}=-\frac{2e^2Y_2}{r^4Y_1}=\frac{1}{\sin^2\theta}\mathcal{Y}^{3}_{24,24};\\
\end{cases}
		\end{eqnarray}	
	
	
\begin{eqnarray}\label{WK}
\begin{cases}
	\mathcal{Y}^{4}_{1223,13}=-\frac{3Y_6^2}{r^6}=\frac{1}{\sin^2\theta}\mathcal{Y}^{4}_{1224,14}=-\mathcal{Y}^{4}_{1213,23}=\frac{1}{\sin^2\theta}\mathcal{Y}^{4}_{1214,24}, \\ 
		\mathcal{Y}^{4}_{1434,13}=\frac{3Y_1Y_6^2\sin^2\theta}{r^6}=-\mathcal{Y}^{4}_{1334,14}, \\
		\mathcal{Y}^{4}_{2434,23}=-\frac{3Y_6^2\sin^2\theta}{r^2Y_1}=-\mathcal{Y}^{4}_{2334,24};
\end{cases}
		\end{eqnarray}

	\begin{eqnarray}\label{KR}
	\begin{cases}
	\mathcal{Y}^{5}_{1223,13}=-\frac{Y_3Y_4}{r^6}=\frac{1}{\sin^2\theta}\mathcal{Y}^{5}_{1224,14}=-\mathcal{Y}^{5}_{1213,23}=\frac{1}{\sin^2\theta}\mathcal{Y}^{5}_{1214,24}, \\ 
		\mathcal{Y}^{5}_{1434,13}=-\frac{Y_1Y_2Y_5\sin^2\theta}{r^6}=-\mathcal{Y}^{5}_{1334,14}, \\
		\mathcal{Y}^{5}_{2434,23}=-\frac{Y_2Y_5\sin^2\theta}{r^2Y_1}=-\mathcal{Y}^{5}_{2334,24};
\end{cases}
		\end{eqnarray}

\begin{eqnarray}\label{KC}
\begin{cases}
	\mathcal{Y}^{6}_{1223,13}=-\frac{Y_2Y_6}{r^6}=\frac{1}{\sin^2\theta}\mathcal{Y}^{6}_{1224,14}=-\mathcal{Y}^{6}_{1213,23}=-\frac{1}{\sin^2\theta}\mathcal{Y}^{6}_{1214,24}, \\ 
	\mathcal{Y}^{6}_{1434,13}=\frac{3Y_2Y_6}{r^6}=-\mathcal{Y}^{6}_{1334,14},\\
	\mathcal{Y}^{6}_{2434,23}=\frac{3Y_2Y_6}{r^2Y_1}=-\mathcal{Y}^{6}_{2334,24}; 
\end{cases}
		\end{eqnarray}	
	
	\begin{eqnarray}\label{KP}
	\begin{cases}
	    \mathcal{Y}^{7}_{1223,13}=-\frac{Y_2Y_4}{r^6}=\frac{1}{\sin^2\theta}\mathcal{Y}^{7}_{1224,14}=-\mathcal{Y}^{7}_{1213,23}=-\frac{1}{\sin^2\theta}\mathcal{Y}^{7}_{1214,24},\\
	    \mathcal{Y}^{7}_{1232,13}=\frac{Y_2 Y_7}{3r^6}=-\frac{1}{\sin^2\theta}\mathcal{Y}^{7}_{1242,14}=\mathcal{Y}^{7}_{1231,23}\\=\mathcal{Y}^{7}_{1321,23}=-\frac{1}{\sin^2\theta}\mathcal{Y}^{7}_{1241,24}=-\frac{1}{\sin^2\theta}\mathcal{Y}^{7}_{1421,24},\\
	    \mathcal{Y}^{7}_{1333,13}=\frac{2e^2 Y_1 Y_2}{3r^6}=\frac{1}{\sin^2\theta}\mathcal{Y}^{7}_{1444,14},\\
	    \mathcal{Y}^{7}_{1434,13}=\frac{Y_1 Y_6 Y_8\sin^2\theta}{3r^6} =-\mathcal{Y}^{7}_{1343,14},\\
	    
	 \mathcal{Y}^{7}_{1443,13}=\frac{Y_1Y_2Y_5\sin^2\theta}{r^6}=-\mathcal{Y}^{7}_{3441,13}=-\mathcal{Y}^{7}_{1334,14},\\
	     \mathcal{Y}^{7}_{2333,23}=-\frac{2e^2Y_2}{3r^2Y_1} =\frac{1}{\sin^2\theta}\mathcal{Y}^{7}_{2444,24},\\
	    \mathcal{Y}^{7}_{2434,23}=\frac{Y_2Y_8\sin^2\theta}{3r^2Y_1}=\mathcal{Y}^{7}_{2343,24}\\
	    \mathcal{Y}^{7}_{2443,23}=\frac{Y_5Y_8\sin^2\theta}{r^2Y_1}=-\mathcal{Y}^{7}_{3442,23} =\mathcal{Y}^{7}_{2334,24};\\
	\end{cases}
			\end{eqnarray}
	
		\begin{eqnarray}\label{KW}
		\begin{cases}
		\mathcal{Y}^{8}_{1223,13}=-\frac{Y_3Y_4}{r^6}=\frac{1}{\sin^2\theta}\mathcal{Y}^{8}_{1224,14}=-\mathcal{Y}^{8}_{1213,23}=\frac{1}{\sin^2\theta}\mathcal{Y}^{8}_{1214,24}, \\ 
			\mathcal{Y}^{8}_{1434,13}=-\frac{Y_1Y_2Y_5\sin^2\theta}{r^6}=-\mathcal{Y}^{8}_{1334,14}, \\
			\mathcal{Y}^{8}_{2434,23}=-\frac{Y_2Y_5\sin^2\theta}{r^2Y_1}=-\mathcal{Y}^{8}_{2334,24};
	\end{cases}
			\end{eqnarray}
	
	\begin{eqnarray}\label{KK}
	\begin{cases}
		\mathcal{Y}^{9}_{1223,13}=-\frac{Y_2Y_6}{r^6}=\frac{1}{\sin^2\theta}\mathcal{Y}^{9}_{1224,14}=-\mathcal{Y}^{9}_{1213,23}=\frac{1}{\sin^2\theta}\mathcal{Y}^{9}_{1214,24}, \\ 
			\mathcal{Y}^{9}_{1434,13}=\frac{3Y_1Y_2Y_6}{r^6}=-\mathcal{Y}^{9}_{1334,14},\\
			\mathcal{Y}^{9}_{2434,23}=\frac{3Y_2Y_6}{r^2Y_1}=-\mathcal{Y}^{9}_{2334,24};
		\end{cases}
				\end{eqnarray}


\begin{eqnarray}\label{Q(g,K)}
\begin{cases}	
	\mathcal{J}^8_{1223,13}=-\frac{3Y_6}{r^2}=\frac{1}{\sin^2\theta}\mathcal{J}^8_{1224,14}=-\mathcal{J}^8_{1213,23}=-\frac{1}{\sin^2\theta}\mathcal{J}^8_{1214,24},\\
	\mathcal{J}^8_{1434,13}=-\frac{3Y_1Y_6\sin^2\theta}{r^2}=-\mathcal{J}^8_{1334,14}, \\
	\mathcal{J}^8_{2434,23}=-\frac{3r^2Y_6\sin^2\theta}{Y_1}=-\mathcal{J}^8_{2334,24}. \\
\end{cases}
		\end{eqnarray}

	
By examining the expressions in equations \eqref{RK}, \eqref{CK}, \eqref{KS}, \eqref{WK}, \eqref{KR}, \eqref{KC}, \eqref{KP}, \eqref{KW}, \eqref{KK} and \eqref{Q(g,K)}, the following results are obtained:

%
%
	 \begin{pr}\label{pr5}
The RNdS spacetime satisfies the following of the pseudosymmetric type conditions:
			\begin{enumerate}[label=(\roman*)]
				\item
		   $\ R\cdot K=-\frac{e^2-Mr+\Lambda r^4}{r^4}Q(g,K),$
		
	\item  $C\cdot K=-\frac{e^2-Mr}{r^4}Q(g,K),$
		
\item  $W\cdot K=-\frac{e^2-Mr}{r^4}Q(g,K),$				
			\item  $K\cdot R=-\frac{e^2-Mr+2\Lambda r^4}{r^4}Q(g,R)\ \mbox{and} \ K\cdot Ric=-\frac{e^2-Mr+2\Lambda r^4}{r^4}Q(g, Ric),$	
\item  $K\cdot C=-\frac{e^2-Mr+2\Lambda r^4}{r^4}Q(g,C)\ \mbox{and} \ K\cdot P=-\frac{e^2-Mr+2\Lambda r^4}{r^4}Q(g, P),$				
			\item  $K\cdot K=-\frac{e^2-Mr+2\Lambda r^4}{r^4}Q(g,K)\ \mbox{and} \ K\cdot W=-\frac{e^2-Mr+2\Lambda r^4}{r^4}Q(g,W).$		
		 	\end{enumerate}	
	\end{pr}

Using Proposition \ref{pr1}--Proposition \ref{pr5}, the following observations regarding the curvature characteristics of the RNdS spacetime can be established:

\begin{thm}
The RNdS spacetime exhibits the following geometric features

	\begin{enumerate}[label=(\roman*)]
		
		\item it is pseudosymmetric as $R\cdot R=-\frac{e^2-Mr+\Lambda r^4}{r^4}Q(g,R).$ \\ Hence $R\cdot C=-\frac{e^2-Mr+\Lambda r^4}{r^4}Q(g,C)$,\ \ $R\cdot P=-\frac{e^2-Mr+\Lambda r^4}{r^4}Q(g,P)$, \\ $R\cdot W=-\frac{e^2-Mr+\Lambda r^4}{r^4}Q(g,W)$\ \  and $R\cdot K=-\frac{e^2-Mr+\Lambda r^4}{r^4}Q(g,K)$;
		\item it is pseudosymmetric due to conformal curvature as  $C\cdot R=-\frac{e^2-Mr}{r^4}Q(g,R).$ \\ Hence $C\cdot C=-\frac{e^2-Mr}{r^4}Q(g,C)$,\ \ $C\cdot P=-\frac{e^2-Mr}{r^4}Q(g,P)$, \\ $C\cdot W=-\frac{e^2-Mr}{r^4}Q(g,W)$\ \  and $C\cdot K=-\frac{e^2-Mr}{r^4}Q(g,K)$;
		\item it admits Ricci pseudosymmetric as $R\cdot Ric=-\frac{e^2-Mr+\Lambda r^4}{r^4}Q(g,Ric).$ \\ Hence $C\cdot Ric=-\frac{e^2-Mr}{r^4}Q(g,Ric)$,\ \ $P\cdot Ric=-\frac{e^2-Mr+\Lambda r^4}{r^4}Q(g, Ric)$, \\ $W\cdot Ric=-\frac{e^2-Mr}{r^4}-\frac{e^2-Mr}{r^4}Q(g,Ric)$\ \  and $K\cdot Ric=-\frac{e^2-Mr+2\Lambda r^4}{r^4}Q(g,Ric)$;

	\item it realizes pseudosymmetric due to concircular curvature as $W\cdot R=-\frac{e^2-Mr}{r^4}Q(g,R).$ \\ Hence $W\cdot C=-\frac{e^2-Mr}{r^4}Q(g,C)$,\ \ $W\cdot P=-\frac{e^2-Mr}{r^4}Q(g,P)$, \\ $W\cdot W=-\frac{e^2-Mr}{r^4}Q(g,W)$\ \  and $W\cdot K=-\frac{e^2-Mr}{r^4}Q(g,K)$;
		 
	\item it exhibits pseudosymmetric due to conharmonic curvature as $K\cdot R=-\frac{e^2-Mr+2\Lambda r^4}{r^4}Q(g,R).$ \\ Hence $K\cdot C=-\frac{e^2-Mr+2\Lambda r^4}{r^4}Q(g,C)$,\ \ $K\cdot P=-\frac{e^2-Mr+2\Lambda r^4}{r^4}Q(g,P)$, \\ $K\cdot W=-\frac{e^2-Mr+2\Lambda r^4}{r^4}Q(g,W)$\ \  and $K\cdot K=-\frac{e^2-Mr+2\Lambda r^4}{r^4}Q(g,K)$;

\item it adheres to the pseudosymmetric type curvature condition  $R\cdot R-Q(Ric,R)=\mathcal{L}_1Q(g,C),$ where $\mathcal{L}_1=\frac{2e^4-6e^2Mr+3M^2r^2-6\Lambda e^2r^4+6M\Lambda r^5)}{3r^4(Mr-e^2)},$
		
	\item the tensor $C\cdot R-R\cdot C$  can be expressed as a linear combination of the tensors $Q(g,C)$, $Q(Ric,C)$, $Q(g,R)$ and $Q(Ric,R)$,

		 

	\item it is an $2$-quasi Einstein manifold, for $L  = -\frac{3\Lambda r^4-e^2}{r^2} $,

		\item it is an Ein$(2),$ as it fulfills $Ric^2+6\Lambda  Ric-\frac{9\Lambda^2r^8-e^4}{r^8} g=0, $

		
		 
		\item the universal form of tensors in RNdS spacetime that are compatible with $R$, $C$, $P$, $W$ and $K$ can be derived as follows: 
		
		$$
		\left(
		\begin{array}{cccc}
			\mathscr{Z}_{11} &\mathscr{Z}_{12} & 0 & 0 \\
			\mathscr{Z}_{12} & \mathscr{Z}_{22} & 0 & 0 \\
			0 & 0 & \mathscr{Z}_{33} & \mathscr{Z}_{34} \\
			0 & 0 & \mathscr{Z}_{34} & \mathscr{Z}_{44}
		\end{array}
		\right),
		$$
		where $\mathscr{Z}_{ij}$ are arbitrary scalars,

	\end{enumerate}
	
\end{thm}

As a specific instance of our finding, we may derive the curvature properties of RNdS spacetimes, which are the generalization of Reissner-Nordstr\"{o}m spacetime, Kottler spacetime, and Schwarzschild spacetime.  The reduction of RNdS spacetime to Reissner-Nordstrom spacetime for $e\neq 0$ and $\Lambda=0$ in \eqref{RNdSM} allows us to propose the following:
\begin{cor}[\cite{Kowa06, SDC}] The Reissner-Nordstr\"{o}m spacetime fulfills the following curvature properties:

\begin{enumerate}[label=(\roman*)]
 \item $\kappa=0$ and consequently $\kappa = C$, $R=W$,
 \item it satisfies $C\cdot R=q_1 Q(g,R)$ and $R\cdot R=q_1 Q(g,R)$, where $q_1=-\frac{e^2-Mr}{r^4}$,
 \item for $q_2=\frac{2e^4-6Me^2r+3M62r^2}{3r^4(Mr^4-e^2)}$ it recognizes $R\cdot R-q_2 Q(g,C)=Q(Ric, R),$
 \item it is Roter type spacetime as $R=\frac{3r}{4e^4}(e^2-Mr) Ric\wedge Ric +\frac{1}{2}g\wedge Ric+\frac{e^2-Mr}{4r^4} g\wedge g,$
 \item it is Ein(2) spacetime as $Ric^2=\frac{e^4}{r^8}g,$
 \item it is $2$-quasi Einstein as rank of $(Ric-L g)$ is $2$ for $L=\frac{e^2}{r^4},$ 
 \item its Ricci tensor is Riemann compatible, projective compatible as well as Weyl compatible.
\end{enumerate}

\end{cor}
If $e=0$ and $\Lambda \neq 0$ in \eqref{RNdSM}, then RNdS spacetime turns into Kottler spacetime(see, \cite{Kowa06}, Exanple $5.2(ii)$).\\
\begin{cor}[\cite{Kowa06}] The curvature properties of the Kottler spacetime are given as follows: 

\begin{enumerate}[label=(\roman*)] 
 \item its scalar curvature is non-zero,
 \item it admits conformal pseudosymmetric condition as $C\cdot C=q_3 Q(g, C),$\\ where $q_3=\frac{3M^2(-8M+5r-8\Lambda r^3)+Mr^2(1+3\Lambda r^2-6\Lambda^2r^4)+r^3(2\Lambda r^2-1)}{6r^3(2M-r+\Lambda r^3)^2},$
 \item it fulfills $R\cdot R-Q(Ric, R)=q_4 Q(g, C),$\\
 where $q_4=\frac{24M^4+2M62r^2(36M\Lambda r+27\Lambda^2 r^4-5)-4M\Lambda r^5(\Lambda^2 r^4-1)-2\Lambda^2 r^4}{3M^2(-8M+5r-8\Lambda r^3)+Mr^2(1+3\Lambda r^2-6\Lambda^2r^4)+r^3(2\Lambda r^2-1)},$
 \item it is a non-Ricci flat Einstein manifold.
\end{enumerate}

\end{cor}

Once again, the RNdS spacetime \eqref{RNdSM} dissolves to the Schwarzschild spacetime if $e=0$ and $\Lambda =0$,  As a result, we have:
\begin{cor}[\cite{SHSB_VBDS}] The curvature properties of the Schwarzschild spacetime are described as follows:

\begin{enumerate}[label=(\roman*)]
 \item it is Ricci flat and as a result $R=W=C=K=P$,
 \item it exceeds $R\cdot R=\frac{M}{r^3} Q(g, R)$ and hence it is a Deszcz pseudosymmetric manifold,
 \item the  curvature of the spacetime is harmonic since it fulfills $divR=0$.
 \item The general form of the compatible tensor for $R$ is as follows: 
 	$$
 		\left(
 		\begin{array}{cccc}
 			\mathscr{Z}_{11} &\mathscr{Z}_{12} & 0 & 0 \\
 			\mathscr{Z}_{12} & \mathscr{Z}_{22} & 0 & 0 \\
 			0 & 0 & \mathscr{Z}_{33} & \mathscr{Z}_{34} \\
 			0 & 0 & \mathscr{Z}_{34} & \mathscr{Z}_{44}
 		\end{array}
 		\right),
 		$$
 		where $\mathscr{Z}_{ij}$ are arbitrary scalars,	
\end{enumerate}

\end{cor}

\begin{rem}
By analyzing various components of different curvature conditions \eqref{RR} to \eqref{Q(g,K)}, we note that RNdS spacetime does not obey any of the following:
	\begin{enumerate}[label=(\roman*)]
		\item since $\nabla R\neq 0$, $\nabla W\neq 0,$ $\nabla P\neq 0$, $\nabla K\neq 0$, $\nabla C\neq 0$,
		


		
		
		
		\item it is not venzi-space for any of the curvature $R, W, P$ and $K,$ 
		
\item 	the Ricci tensor of the RNdS spacetime neither cyclic parallel nor Codazzi-type,	
		\item the curvature $2$-forms associated with $R$, $W,$ $P$ and $K$ are not recurrent. 
		
		\item neither locally symmetric nor semisymmetric,
		\item the RNdS spacetime does not exhibit Chaki pseudosymmetry or weak symmetry for $R$ $	W,$ $P,$ $C,$ and $K$. 
		
	\end{enumerate}  
\end{rem}

\section{\bf Nature of the Energy momentum tensor of the Reissner-Nordstr\"{o}m-de Sitter spacetime}

In the framework of general relativity, Albert Einstein described the structure of spacetime through a set of field equations that capture its geometric behavior. The fundamental equation takes the form
$$Ric - \frac{\kappa}{2}g + \Lambda g = \frac{8\pi G}{c^4}T,$$
where $Ric$ denotes the Ricci curvature tensor, $\kappa$ is the scalar curvature, and $T$ stands for the energy-momentum tensor of the spacetime. Here, $\Lambda$ represents the cosmological constant, $G$ is the gravitational constant, and $c$ is the speed of light in vacuum. This relation mathematically expresses how energy and matter influence the curvature of spacetime, playing a central role in our understanding of gravity and the dynamics of the universe.



Assuming $\frac{8\pi G}{c^4} = 1$, the components of the energy-momentum tensor can be written as:

$$\begin{array}{c}
	T_{11}=\frac{(3\Lambda r^4+e^2) Y_1}{r^6}, \ \
	T_{22}=-\frac{3\Lambda r^4+e^2}{r^2Y_1}, \\
	T_{33}=\frac{3\Lambda r^4-e^2}{r^2}=\frac{1}{\sin^2\theta}T_{44}. \\	
\end{array}$$
The non-zero components of the covariant derivatives of $T$ are given below:
\begin{eqnarray*}
 \begin{cases}
 T_{11,2}=-\frac{4e^2Y_1}{r^7},\\
 T_{22,2}=\frac{4e^2}{r^3Y_1},\\
 T_{23,3}=\frac{2e^2}{r^3}=\frac{1}{\sin^2\theta}T_{24,4}'\\
 T_{33,2}=\frac{4e^2}{r^3}=\frac{1}{\sin^2\theta}T_{44,2}.\\
 \end{cases}
\end{eqnarray*} 

The following are the non-vanishing components of $R\cdot T$, $C\cdot T$, $W\cdot T$, and $K\cdot T$:
\begin{eqnarray}
 \begin{cases}
 (R\cdot T)_{1313}=\frac{2e^2Y_1Y_2}{r^8}=\frac{1}{\sin^2\theta}(R\cdot T)_{1414},\\
 (R\cdot T)_{2323}=-\frac{2e^2Y_2}{r^4Y_1}=\frac{1}{\sin^2\theta}(R\cdot T)_{2424}.\\
 \end{cases}
\end{eqnarray}

\begin{eqnarray}
 \begin{cases}
 (C\cdot T)_{1313}=\frac{2e^2Y_1Y_6}{r^4}=\frac{1}{\sin^2\theta}(C\cdot T)_{1414},\\
 (C\cdot T)_{2323}=\frac{2e^2Y_6}{r^4Y_1}=\frac{1}{\sin^2\theta}(C\cdot T)_{2424}.\\
 \end{cases}
\end{eqnarray}

\begin{eqnarray}
 \begin{cases}
 (W\cdot T)_{1313}=\frac{2e^2Y_1Y_6}{r^4}=\frac{1}{\sin^2\theta}(W\cdot T)_{1414},\\
 (W\cdot T)_{2323}=\frac{2e^2Y_6}{r^4Y_1}=\frac{1}{\sin^2\theta}(W\cdot T)_{2424}.\\
 \end{cases}
\end{eqnarray}

\begin{eqnarray}
 \begin{cases}
 (K\cdot T)_{1313}=\frac{2e^2Y_1Y_2}{r^8}=\frac{1}{\sin^2\theta}(K\cdot T)_{1414},\\
 (K\cdot T)_{2323}=-\frac{2e^2Y_2}{r^4Y_1}=\frac{1}{\sin^2\theta}(K\cdot T)_{2424}.\\
 \end{cases}
\end{eqnarray}
The non-vanishing components of the tensor $Q(g,T)$ are given below. 
\begin{eqnarray}
 \begin{cases}
 Q(g, T)_{1313}=-\frac{2e^2Y_1}{r^4}=\frac{1}{\sin^2\theta}Q(g, T)_{1414},\\
 Q(g, T)_{2323}=\frac{2e^2}{Y_1}=\frac{1}{\sin^2\theta}Q(g, T)_{2424}.\\
 \end{cases}
\end{eqnarray}

\indent By virtue of the components listed above, the following conclusions can be drawn:
\begin{thm} The characteristics of the energy-momentum tensor in the RNdS spacetime are examined as follows:

\begin{enumerate}[label=(\roman*)]
\item  $R\cdot  T=-\frac{e^2-Mr+\Lambda r^4}{r^4} Q(g,T)$, 
		\item $C\cdot  T=-\frac{e^2-Mr}{r^4} Q(g,T)$,
		
		\item $W\cdot  T=-\frac{e^2-Mr}{r^4} Q(g,T)$,
			\item $K\cdot  T=-\frac{e^2-Mr+2\Lambda r^4}{r^4} Q(g,T)$,
%
%
%
%
%
%
		\item The energy-momentum tensor is compatible with the structures of projective, conharmonic, concircular, Riemann and conformal curvature.
	\end{enumerate}
\end{thm}

\begin{cor}[] The nature of the energy-momentum tensor of Reissner-Nordstr\"{o}m, Kottler, and Schwarzchild spacetimes are respectively described as follows:

\begin{enumerate}[label=(\roman*)]
 \item Energy momentum tensor of Reissner-Nordstr\"{o}m spacetime:\\
 $(a)$ trace of energy momentum tensor vanishes,\\
 $(b)$ it is neither codazzi type nor cyclic parallel,\\
 $(c)$ it shows some pseudosymmetric condition such as $R\cdot T=\frac{Mr-e^2}{r^4}Q(g,T),$ $C\cdot T=\frac{Mr-e^2}{r^4}Q(g,T),$ $W\cdot T=\frac{Mr-e^2}{r^4}Q(g,T)$ and $K\cdot T=\frac{Mr-e^2}{r^4}Q(g,T),$
 \item Energy momentum tensor of Kottler spacetime:\\
  $(a)$ trace of energy momentum tensor is $12\Lambda$,\\
  $(b)$ it is codazzi type as well as cyclic parallel,\\
  $(c)$ it does not reveal pseudosymmetric condition,
 \item Energy momentum tensor of Schwarzchild spacetime:\\
   $(a)$ trace of energy momentum tensor is zero,\\
   $(b)$ it is codazzi type as well as cyclic parallel,\\
   $(c)$ it does not fulfill pseudosymmetric condition.
 \end{enumerate}

\end{cor}

%

%

\section{\bf Ricci soliton and symmetries on Reissner-Nordstr\"{o}m-de Sitter spacetime}


 Within the RNdS spacetime, the vector fields $\frac{\partial}{\partial t}$ and $\frac{\partial}{\partial \phi}$ are identified as Killing vector fields, indicating that the Lie derivatives satisfy $\pounds_{\partial/\partial t} g = 0$ and $\pounds_{\partial/\partial \phi} g = 0,$ where $\pounds$ represents the Lie derivative operator. We note that the coordinate vector fields ${\frac{\partial}{\partial r}}$ and  $\frac{\partial}{\partial \theta}$ are not Killing. 

 Defining $\mathcal{A} = \pounds_{\frac{\partial}{\partial r}} g$, the corresponding non-zero components of $\mathcal{A}$ are computed as follows:
$$\begin{array}{c}
	\mathcal A_{11}=\frac{2Y_2}{r^3},\ \ \ 
	\mathcal A_{22}=\frac{2rY_2}{Y_1^2},
	\nonumber\\
	\mathcal A_{33}=2r,\ \ \ 
	\mathcal A_{44}=2r\sin^2\theta.\nonumber
\end{array}$$

Therefore, in the RNdS spacetime, for the non-Killing vector field $\frac{\partial}{\partial r}$ and the 1-form $\eta = (0,1,0,0)$, the following relation holds:

\begin{eqnarray}\nonumber
	\pounds_{\frac{\partial}{\partial r}}g+ 2\sigma_1 Ric+2\sigma_2 g-2\sigma_3 \eta\otimes\eta=0,
\end{eqnarray}
where $\sigma_1,\sigma_2,\sigma_3$ are given by
\begin{equation}\label{LieCoefficient}
	\left.
	\begin{aligned}
		\sigma_1&=\frac{r^3(3Mr-r^2-2e^2)}{2e^2Y_1},\\
		\sigma_2&=\frac{-9M\Lambda r^4+3\Lambda r^5+Me^2-re^2+8\Lambda r^3e^2}{2e^2Y_1},\\
		\sigma_3&=-\frac{2(Mr^2-\Lambda r^5-re^2)}{Y_1^2}.
	\end{aligned}\ \ 
	\right\rbrace	
\end{equation}
This leads to the following:

\begin{thm} \label{thm6.1}
	In the RNdS spacetime if $(2Mr - r^2 + \Lambda r^4 - e^2)$ is non-zero then it reveals an almost $\eta$-Ricci-Yamabe soliton for the soliton vector field $\xi=\frac{\partial}{\partial r}$ and for the $1$-form $\eta = (0,1,0,0)$ satisfying:
	$$\frac{1}{2}\pounds_\xi g+\sigma_1 Ric+\left(\lambda-\frac{1}{2}\sigma_4\kappa\right) g+\sigma_3 \eta\otimes\eta=0,$$
	where $\sigma_4=2$,   $\lambda=\sigma_2+\kappa$, and $\sigma_1,\sigma_2,\sigma_3$ are given in (\ref{LieCoefficient}).
	
\end{thm}

\begin{thm}\label{thm6.2}
If $r^3(3Mr - r^2 - 2e^2) = 2Mr - r^2 + \Lambda r^4 - e^2$ and $(2Mr - r^2 + \Lambda r^4 - e^2) \neq 0$, then the RNdS spacetime fulfills an almost $\eta$-Ricci soliton for the soliton vector field $\xi=\frac{\partial}{\partial r}$ and the 1-form $\eta = (0,1,0,0)$ satisfying:
	\begin{eqnarray}
		&&\frac{1}{2}\pounds_{\xi}g+Ric+\sigma_2 g-\sigma_3\eta\otimes\eta=0\nonumber,
	\end{eqnarray}
where $\sigma_2$, $\sigma_3$ is given in (\ref{LieCoefficient}).
\end{thm}

\begin{rem}
It is observed that if the RNdS spacetime fulfills the condition $Mr^2 - \Lambda r^5 - r e^2 = 0$ and $r^4+r^3(3M-2e^2)+2e^2r^2(\Lambda r^3-1)+4 Me^2 r-2e^4 = 0$, then, based on Theorem \eqref{thm6.1} and Theorem \eqref{thm6.2}, the spacetime realizes a Ricci soliton with regard the soliton vector field $\xi = \frac{\partial}{\partial r}$.
\end{rem}

Define $\mathcal{M}^1 = \pounds_{\frac{\partial}{\partial r}} Ric$, $\mathcal{G} = \pounds_{\frac{\partial}{\partial r}} \widetilde{R}$, and $\mathcal{M}^2 = \pounds_{\frac{\partial}{\partial r}} R$. The corresponding non-zero components of $\mathcal{M}^1$, $\mathcal{G}$, and $\mathcal{M}^2$ are calculated as follows:
$$\begin{array}{c}
 \mathcal{M}^1_{11}=-\frac{\{6\Lambda r^5(M-\Lambda r^3)-2r(5M-2r+4\Lambda r^3)e^2+6e^2 \} }{r^7},\\
\mathcal{M}^1_{22}=-\frac{\{6\Lambda r^5(-M+\Lambda r^3)+3r(-3M+2r)e^2+3e^4 \} }{r^3Y_1^2},\\
	\mathcal{M}^1_{33}=\frac{\{2(3\Lambda+e^2)\}}{r^3}=\frac{1}{\sin^2\theta} \mathcal{M}^1_{44};
\end{array}$$

$$\begin{array}{c}
	\mathcal G^1_{212}=\frac{2r\left\lbrace 4M^2+\Lambda^2r^6+Mr(-3+4\Lambda r^2)\right\rbrace-4r(5M-3r+4\Lambda r^3)e^2+6e^4}{r^3Y_1^2}=-\mathcal G^1_{221}, \\ 
	\mathcal G^1_{313}=\frac{Mr+2\Lambda r^4-2e^2}{r^3}=-\mathcal G^1_{331}=\mathcal G^2_{323}=-\mathcal G^1_{332} \\
	=\frac{1}{\sin^2\theta}\mathcal G^1_{414}=-\frac{1}{\sin^2\theta}\mathcal G^1_{441}=\frac{1}{\sin^2\theta}\mathcal G^1_{424}=-\frac{1}{\sin^2\theta}\mathcal G^1_{442},\\
	
	\mathcal G^2_{112}=\frac{2r\left\lbrace 8M^2+\Lambda^2r^6-Mr(-3-2\Lambda r^2)\right\rbrace-4r(10M-3r+2\Lambda r^3)e^2+18e^4}{r^7}=-\mathcal G^2_{121}, \\
	
	\mathcal G^3_{113}=-\frac{-8M^2r^2+M(3r^2+\Lambda r^5+15re^2)-2(\Lambda^2 r^8+2r^2e^2+3e^4)}{r^7}=-\mathcal G^3_{131}=\mathcal G^4_{114}=-\mathcal G^4_{141}, \\
	
	\mathcal G^3_{223}=\frac{r^2\left\lbrace 4M^2-2\Lambda^2r^6+Mr(-3+7\Lambda r^2)\right\rbrace+r^2(-7M+4r-8\Lambda r^3)e^2+2e^4}{r^3Y_1^2}=\mathcal G^4_{224}=-\mathcal G^3_{232}=-\mathcal G^4_{242},\\

	\mathcal G^3_{434}=\frac{2Y_2\sin^2\theta}{r^3}=-\mathcal G^3_{443}=-\sin^2\theta\mathcal G^4_{334}=\sin^2\theta\mathcal G^4_{343};\\ 
\end{array}$$

$$\begin{array}{c}
	\mathcal{M}^2_{1212}=\frac{6(Mr-2e^2)}{r^5}=-\mathcal{M}^2_{1221}=\mathcal{M}^2_{2121}=-\mathcal{M}^2_{2112}, \\
	\mathcal{M}^2_{1313}=\frac{r^2\left\lbrace 4M^2+Mr(-1+\Lambda r^2)+2\Lambda r^4(-1+2\Lambda r)\right\rbrace+r(-9M+2r)e^2+4e^4}{r^5}=-\mathcal{M}^2_{1331}\\=-\frac{1}{\sin^2\theta} \mathcal{M}^2_{1414}=\frac{1}{\sin^2\theta} \mathcal{M}^2_{1441}=-\mathcal{M}^2_{3113}=\mathcal{M}^2_{3131} \\ =\frac{1}{\sin^2\theta} \mathcal{M}^2_{4141}=-\frac{1}{\sin^2\theta} \mathcal{M}^2_{4114}, \\
	\mathcal{M}^2_{2323}=\frac{r^2\left\lbrace M-9M\Lambda r^2+2\Lambda r^3 \right\rbrace+\left\lbrace M-2r+8\Lambda r^3\right\rbrace e^2}{Y_1^2 }=-\mathcal{M}^2_{2332}=\frac{1}{\sin^2\theta} \mathcal{M}^2_{2424} \\=-\frac{1}{\sin^2\theta} \mathcal{M}^2_{2442}=-\mathcal{M}^2_{3223}= \mathcal{M}^2_{3232}=\frac{1}{\sin^2\theta} \mathcal{M}^2_{4242}=-\frac{1}{\sin^2\theta} \mathcal{M}^2_{4224}, \\
	\mathcal{M}^2_{3434}=2(M+2\Lambda r^3)\sin^2\theta=-\mathcal{M}^2_{3443}=-\mathcal{M}^2_{4334}=\mathcal{M}^2_{4343}.
\end{array}$$

Let $\mathcal{M}^3 = \pounds_{\frac{\partial}{\partial \theta}} g$, $\mathcal{M}^4 = \pounds_{\frac{\partial}{\partial \theta}} Ric$, $\mathcal{Q} = \pounds_{\frac{\partial}{\partial \theta}} \widetilde{R}$, and $\mathcal{M}^5 = \pounds_{\frac{\partial}{\partial \theta}} R$. Then, the non-vanishing components of $\mathcal{M}^3$, $\mathcal{M}^4$, $\mathcal{Q}$, and $\mathcal{M}^5$ are as follows:
$$\begin{array}{c}
\mathcal{M}^3_{11}=0, \ \ \mathcal{M}^3_{22}=0, \ \ \mathcal{M}^3_{33}=0, \ \ \mathcal{M}^3_{44}=r^2 \sin2\theta; \\
\end{array}$$
$$\begin{array}{c}
\mathcal{M}^4_{11}=0, \ \ \mathcal{M}^4_{22}=0, \ \ \mathcal{M}^4_{33}=0, \ \ \mathcal{M}^4_{44}=\frac{(3\Lambda r^4+e^2) \sin2\theta}{r^2}; \\
\end{array}$$
$$\begin{array}{c}
	\mathcal Q^1_{414}=\frac{Y_2\sin2\theta}{r^2}=-\mathcal Q^1_{441}=\mathcal Q^2_{424}=-\mathcal Q^2_{442}, \ \
	\mathcal Q^3_{434}=\frac{(2Mr+\Lambda r^4-e^2)\sin2\theta}{r^2}=-\mathcal Q^3_{443};
\end{array}$$

$$\begin{array}{c}
		\mathcal{M}^5_{1414}=\frac{Y_1Y_2 \sin2\theta}{r^4}=-\mathcal{M}^5_{1441}=-\frac{1}{\sin^2\theta}\mathcal{M}^5_{4114}=\frac{1}{\sin^2\theta}\mathcal{M}^5_{4141}, \\
	\mathcal{M}^5_{2424}=\frac{Y_2\sin2\theta}{Y_1}=-\mathcal{M}^5_{2442}=\mathcal{M}^5_{4224}=\mathcal{M}^5_{4242}, \\
	\mathcal{M}^5_{3434}=(2Mr+\Lambda r^4-e^2)\sin2\theta=-\mathcal{M}^5_{3443}=\mathcal{M}^5_{4334}=\mathcal{M}^5_{4343}.
\end{array}$$

Based on the preceding computations of the Lie derivatives of various curvature tensors, it is evident that the RNdS spacetime does not reveal any of the following symmetry due to the vector fields $\frac{\partial}{\partial r}$, $\frac{\partial}{\partial \theta}$,
	 \begin{enumerate}[label=(\roman*)]
	 	\item it neither admits Ricci collineation nor Ricci inheritance,
	 	
	 	\item it does not possesses curvature collineation with regard to curvature tensor of type $(1,3)$ and $(0,4)$,
	 	
	 	\item it also does not support curvature inheritance for either the $(1,3)$-type or the $(0,4)$-type curvature tensors.
	 \end{enumerate}

\section{\bf Reissner-Nordstr\"{o}m-de Sitter spacetime Vs Vaidya-Bonner-de Sitter spacetime  }
The Vaidya-Bonner-de Sitter(briefly,VBdS) spacetime (\cite{VBdS_1987}) serves as an exact solution to Einstein's field equations incorporating a non-zero cosmological constant. It describes the gravitational field generated by a spherically symmetric mass distribution. In advanced Eddington–Finkelstein coordinates $(t, r, \theta, \phi)$, the line element for the VBdS spacetime is expressed as:
\begin{eqnarray}\label{VBdS}
	ds^2=\left(1-\frac{2m(t)}{r}+\frac{q^2(t)}{r^2}-\frac{\lambda r^2}{3}\right)dt^2-2dtdr-r^2(d\theta^2+\sin^2\theta d\phi^2),
    \end{eqnarray}
     where $\lambda$ is the cosmological constant, $t$ is the advanced time coordinate, $r$ is the radial coordinate, and both the mass $m(t)$ and charge $q(t)$ of the body depend on time. The detailed study of VBdS spacetime are discussed in \cite{SHSB_VBDS}. A comparison of the curvature properties between RNdS spacetime and VBdS spacetime is outlined as follows:\\

\noindent\textbf{Similarities:}
\begin{enumerate}[label=(\roman*)]
\item scalar curvature is non-zero for both spacetimes,
	\item both spacetimes describes $2$-quasi-Einstein manifolds,
	
	
	\item the  coordinate  vector fields $\frac{\partial}{\partial \phi}$ is Killing for both spacetimes,
	
	\item the vector fields $\frac{\partial}{\partial r}$ and $\frac{\partial}{\partial \theta}$ are non-Killing vectors in both spacetimes.
\item the Ricci tensors of both spacetimes are neither of the Codazzi type nor cyclic parallel,
\item  the Ricci tensor exhibits compatibility with the Riemann, conformal, concircular, conharmonic, and Weyl projective curvature tensors.
	
\item Under certain conditions both spacetimes admits almost  Ricci soliton for the non-Killing vector field $\frac{\partial}{\partial r},$
	\item the energy momentum tensors of both spacetimes are Riemann compatible, projective compatible, conharmonic compatible, concircular compatible and conformal compatible.
\end{enumerate}
Nevertheless, they exhibit differences in the following aspects:\\

\noindent\textbf{Dissimilarities:}

\begin{enumerate}[label=(\roman*)]
	\item the conformal $2$-forms of RNdS spacetime are recurrent while VBdS spacetime does not admit such recurrence,
	
	
	
		
	\item The energy-momentum tensor in the RNdS spacetime exhibits pseudosymmetry, including pseudosymmetry with respect to the conformal, concircular, and conharmonic curvature tensors. In contrast, the energy-momentum tensor in the VBdS spacetime does not admit such pseudosymmetries but satisfies the relation $Q(T, R) = -2\lambda\, Q(g, R) + Q(Ric, R)$ ,
\item the VBdS spacetime realizes generalized Roter type condition, while the RNdS spacetime is of Roter type,
	
	\item in RNdS spacetime the coordinate vector field $\frac{\partial}{\partial t}$ is Killing whereas VBdS spacetime $\frac{\partial}{\partial t}$ is non-Killing,
	
\item The VBdS spacetime fulfills almost $\eta$-Ricci Yamabe soliton for the non-Killing vector field $\frac{\partial}{\partial \theta}$ whereas RNdS spacetimes realizes such notion for the non-Killing vector field $\frac{\partial}{\partial r},$	
	\item The VBdS spacetime obeys generalized conharmonic curvature inheritance corresponding to the non-Killing coordinate vector field $\frac{\partial}{\partial \theta}$, whereas the RNdS spacetime does not satisfy this property for the same non-Killing vector field.

\end{enumerate}
\section{\bf  Acknowledgment}


The Second author greatly acknowledges to The University Grants Commission, Government of India for the award of junior Research Fellow. All the algebraic computations of Section $3-5$ are performed by a program in Wolfram Mathematica developed by the first author A. A. Shaikh.

\section{\bf Declarations}

{\bf Data availability:} Data sharing not applicable to this article as no data sets were used/generated or analyzed during the current study.\\

{\bf Competing interests:} The authors have no competing interests to declare that are relevant to the content of this article.\\

{\bf Funding:} No funding was received to assist with the preparation of this manuscript.\\


%



\begin{thebibliography}{99}\baselineskip=16pt

\bibitem{AD83}
Adam\'{o}w, A. and Deszcz, R.,
\emph{On totally umbilical submanifolds of some class of Riemannian manifolds},
Demonstratio Math.,
\textbf{16} (1983), 39--59.


\bibitem{Ahsan1977_231}
Ahsan, Z.,
\emph{Algebraic classification of space-matter tensor in general relativity},
Indian J. Pure Appl. Math.,
\textbf{8(2)} (1977), 231--237


\bibitem{Ahsan1977_1055}
Ahsan, Z.,
\emph{Algebra of space-matter tensor in general relativity},
Indian J. Pure Appl. Math.,
\textbf{8(9)} (1977), 1055--1061.



\bibitem{Ahsan1978}
Ahsan, Z., \emph{Collineation in electromagnetic field in general relativity- The null field case},
Tamkang J. Maths.,
\textbf{9(2)} (1978), 237.


\bibitem{Ahsan1987} 
Ahsan, Z.,
\emph{On the Nijenhuis tensor for null electromagnetic field},
J. Math. Phys. Sci.,
\textbf{21(5)} (1987), 515--526.


\bibitem{Ahsan1995}
Ahsan, Z., \emph{Symmetries of the electromagnetic fields in general relativity},
Acta Phys. Sinica,
\textbf{4} (1995), 337.


\bibitem{Ahsan1996}
Ahsan, Z.,
\emph{A symmetry property of the space-time of general relativity in terms of the
	space-matter tensor},
Braz. J. Phys. 
\textbf{26(3)} (1996), 572-576.



\bibitem{Ahasan2005}
Ahsan, Z.,
\emph{On a geometrical symmetry of the space-time of general relativity},
Bull. Call. Math. Soc.,
\textbf{97(3)} (2005), 191--200.

\bibitem{Ahsan2018}
Ahsan, Z.,
\emph{Ricci solitons and the spacetime of general relativity},
J. Tensor Soc.,
\textbf{12} (2018), 49--64.

\bibitem{AA2012}
Ahsan, Z. and Ali, M.,
\emph{Symmetries of type D pure radiation fields}.
Int. J. Theo. Phys.,
\textbf{51} (2012), 2044-2055.

\bibitem{AhsanAli2014}
Ahsan, Z. and Ali, M.,
\emph{On some properties of $W$-curvature tensor},
Palestine J. Math.,
\textbf{3(1)} (2014), 61--69.

\bibitem{AH1980}
Ahsan, Z. and Husain, S. I.,
\emph{Null electromagnetic fields, total gravitational radiation
	and collineations in general relativity},
Annali di Mathematical Pura ed Applicata,
\textbf{126} (1980), 379396.

\bibitem{AliAhsan2012}
Ali, M. and Ahsan Z.,
\emph{Ricci solitons and symmetries of spacetime manifold of general relativity},
Global J. Adv. Research Classical Mod. Geom., 
\textbf{1(2)} (2012), 75--84.

\bibitem{AliAhsan2013}
Ali, M. and  Ahsan, Z.,
\emph{Geometry of Schwarzschild soliton},
J. Tensor Soc.,
\textbf{7} (2013), 49--57.

\bibitem{AliAhsan2015}
Ali, M. and Ahsan, Z.,
\emph{Gravitational field of Schwarzschild soliton},
Arab J. Math. Sci.,
\textbf{21(1)} (2015), 15--21.

\bibitem{HAli2003}
Ali, M. H., \emph{Spinning particles in Reissner--Nordstr\"om--de Sitter spacetime}, Gen. Relativity Gravitation \textbf{35(2)}, (2003), 285--305.

\bibitem{ARS95}
Al\'ias, L. J., Romero, A. and S\'anchez, M., \emph{Uniqueness of complete spacelike hypersurfaces of constant mean curvature in generalized Robertson-Walker space-times}, Gen. Relativity Gravitation, \textbf{27(1)} (1995), 71--84.



\bibitem{ADEHM14}
Arslan, K., Deszcz, R., Ezenta\c{s}, R., Hotlo\'{s}, M. and Murathan, C., \emph{On generalized Robertson-Walker spacetimes satisfying some curvature condition}, Turkish J. Math., \textbf{38(2)} (2014), 353--373.




\bibitem{Bess87}
Besse, A. L.,
\emph{Einstein manifolds},
Springer-Verlag, Berlin, Heidelberg, \textbf{1987}.

\bibitem{Blaga2016}
Blaga, A. M.,
\emph{$\eta$-Ricci solitons on Lorentzian para-Sasakian manifolds},
Filomat, 
\textbf{30(2)} (2016), 489--496.
\bibitem{Brink1925}
Brinkmann, H. W., \emph{Einstein spaces which are mapped conformally on each other}, Math. Ann. \textbf{94}
(1925), 119--145.












\bibitem{Cart26}
Cartan, \'E., \emph{Sur une classe remarquable d'espaces de Riemannian}, Bull. Soc. Math. France, \textbf{54} (1926), 214- 264.

\bibitem{Cart46}
Cartan, \'E.,
\emph{Le\c cons sur la g\' eom\' etrie des espaces de Riemann}, 2nd ed., Paris, \textbf{1946}.

\bibitem{Chak87}
Chaki, M. C.,
\emph{On pseudosymmetric manifolds},
An. \c{S}tiin\c{t}.  Univ. AL. I. Cuza Ia\c{s}i. Mat. (N.S.)  Sect. Ia,
\textbf{33(1)} (1987), 53--58.

\bibitem{Chak88}
Chaki, M. C., \emph{On pseudo Ricci symmetric manifolds}, Bulgarian J. Phys., {\bf 15} (1988), 526--531.

\bibitem{C01}
Chaki, M. C., \emph{On generalized quasi-Einstein manifolds,} Publ. Math. Debrecen, \textbf{58} (2001), 683--691.




\bibitem{Cho2009}
Cho, J. and Kimura, M.,
\emph{Ricci solitons and real hypersurfaces in a complex space form}, 
Tohoku Math. J.,
\textbf{61(2)} (2009), 205--212.

\bibitem{CCDNP2023}
Chrysostomou, A., Cornell, A. S., Deandrea, A., Noshad, H., and Park, S. C., \emph{Reissner-Nordstr\"om black holes in de Sitter space-time: Bounds with quasinormal frequencies}, (2023), https://arxiv.org/abs/2310.07311.

\bibitem{CI2025}
Czinner, V. G., Iguchi, H., \emph{Hawking–Rényi black hole thermodynamics, Kiselev solution, and cosmic censorship}, Eur. Phys. J. C \textbf{85}, 443 (2025). https://doi.org/10.1140/epjc/s10052-025-14117-w


\bibitem{DDHKS00}
Defever, F.,  Deszcz, R.,  Hotlo$\acute{\mbox{s}}$, M.,  Kucharski, M.  and  Sent$\ddot{\mbox{u}}$rk, Z.,  \emph{Generalisations of Robertson-Walker spaces}, Ann. Univ. Sci. Budapest, E$\ddot{\mbox{o}}$tv$\ddot{\mbox{o}}$s Sect. Math., {\bf 43} (2000), 13--24.
\bibitem{DD91}
Defever, F. and  Deszcz, R., \emph{On semi-Riemannian manifolds satisfying the condition $R\cdot R=Q(S,R)$}, in:Geometry and Topology of Submanifolds III, World Sci., River Edge, NJ, (1991), 108-130.

\bibitem{DVV1991}
Deszcz, R., Verstrae1en, L. and Vrancken, L., \emph{On the symmetry of warped product spacetimes}, Gen. Relativity Gravitation \textbf{23} (1991), 671--681.
\bibitem{Desz92}
Deszcz, R., \emph{On pseudosymmetric spaces},
Bull. Belg. Math. Soc., Ser. A,
\textbf{44} (1992), 1--34.

\bibitem{Desz93}
Deszcz, R., \emph{Curvature properties of a pseudosymmetric manifolds}, Colloq. Math., \textbf{62} (1993), 139--147.

\bibitem{Desz03}
Deszcz, R., \emph{On Roter type manifolds}, 5-th Conference on Geometry and Topology of Manifolds, Krynica, Poland, April 27 - May 3, (2003),  25.

\bibitem{DG02}
Deszcz, R. and G\l ogowska, M.,
\emph{Some examples of nonsemisymmetric Ricci-semisymmetric hypersurfaces},
Colloq. Math., \textbf{94} (2002), 87--101.

\bibitem{DGHS98}
Deszcz, R., G\l ogowska, M., Hotlo\'{s}, M. and \d Sent\"{u}rk, Z.,
\emph{On certain quasi-Einstein semi-symmetric hypersurfaces},
Ann. Univ. Sci. Budapest E\"{o}tv\"{o}s Sect. Math.,
\textbf{41} (1998), 151--164.

\bibitem{DGJPZ13}
Deszcz, R., G\l ogowska, M., Je\l owicki, L., Petrovi\'{c}-Torga\u{s}ev, M. and Zafindratafa, G.,
\emph{On Riemann and Weyl compatible tensors},
Publ. Inst. Math. (Beograd) (N.S.),
\textbf{94(108)} (2013), 111--124.

\bibitem{DGHS11}
Deszcz, R., G\l ogowska, M., Hotlo\'s, M. and Sawicz, K.,
\emph{A survey on generalized Einstein metric conditions},
Advances in Lorentzian Geometry, Proceedings of the Lorentzian Geometry Conference in Berlin, AMS/IP Studies in Advanced Mathematics, \textbf{49}, S.-T. Yau (series ed.),
M. Plaue, A.D. Rendall and M. Scherfner (eds.), 2011, 27-46.

\bibitem{DGHZ15} 
Deszcz, R., G\l ogowska, M., Hotlo\'s, M. and Zafindratafa, G.,
\emph{On some curvature conditions of pseudosymmetric type}, 
Period. Math. Hungarica,
\textbf{70(2)} (2015), 153--170.

\bibitem{DGHZ16}
Deszcz, R., G\l ogowska, M., Hotlo\'s, M. and Zafindratafa, G., \emph{Hypersurfaces in space
forms satisfying some curvature conditions}, J. Geom. Phys., \textbf{99} (2016), 218--231.

\bibitem{DGJZ-2016}
Deszcz, R., G\l ogowska, M., Je\l owicki, J. and Zafindratafa, Z.,
\emph{Curvature properties of some class of warped product manifolds},
Int. J. Geom. Methods Mod. Phys.,
\textbf{13} (2016), 1550135.

\bibitem{DGP-TV-2015}
Deszcz, R., G\l ogowska, M., Petrovi\'{c}-Torga\u{s}ev, M. and Verstraelen, L.,
\emph{Curvature properties of some class of minimal hypersurfaces in Euclidean spaces},
Filomat,
\textbf{29} (2015), 479--492.

\bibitem{DGPSS11}
Deszcz, R., G\l ogowska, M., Plaue, M., Sawicz, K. and Scherfner, M.,
\emph{On hypersurfaces in space forms satisfying particular curvature conditions of Tachibana type},
Kragujevac J. Math.,
\textbf{35} (2011), 223--247.

\bibitem{DGP-TV-2011}
Deszcz, R., G\l ogowska, M., Petrovi\'c-Torga\u{s}ev, M. and Verstraelen, L.,
\emph{On the Roter type of Chen ideal submanifolds},
Results Math.,
\textbf{59} (2011), 401--413.
\bibitem{DHV2004}	
Deszcz, R., Haesen, S. and Verstraelen, L., \emph{Classification of space-times satisfying some pseudo-
	symmetry type conditions}, Soochow J. Math. \textbf{30} (2004), 339--349.


\bibitem{DH03}
Deszcz, R. and Hotlo\'{s}, M.,
\emph{On hypersurfaces with type number two in spaces of constant curvature},
Ann. Univ. Sci. Budapest E\"{o}tv\"{o}s Sect. Math.,
\textbf{46} (2003), 19--34.

\bibitem{DHJKS14}
Deszcz, R., Hotlo\'{s}, M., Je\l owicki, J., Kundu, H. and Shaikh, A. A.,
\emph{Curvature properties of G\"{o}del metric},
Int. J. Geom. Methods Mod. Phys., \textbf{11} (2014), 1450025. Erratum: \emph{Curvature properties of Gödel metric}, Int. J. Geom. Methods Mod. Phys., \textbf{16} (2019), 1992002.

\bibitem{DK99}
Deszcz, R. and Kucharski, M., \emph{On curvature properties of certain generalized Robertson-Walker spacetimes}, Tsukuba J. Math., \textbf{23(1)} (1999), 113--130.

\bibitem{DPSch-2013}
Deszcz, R., Plaue, M. and Scherfner, M.,
\emph{On Roter type warped products with 1-dimensional fibres},
J. Geom. Phys., \textbf{69} (2013), 1--11.


\bibitem{Duggal1992}
Duggal, K. L., 
\emph{Curvature inheritance symmetry in Riemannian spaces with applications to fluid space times},
J. Math. Phys., \textbf{33(9)} (1992), 2989--2997.

\bibitem{EC21}
Eyasmin, S. and Chakraborty, D., \emph{Curvature properties of (t-z)-type plane wave metric}, J. Geom. Phys., \textbf{160} (2021), 104004.

\bibitem{ECS22}
Eyasmin, S., Chakraborty, D. and Sarkar, M., \emph{Curvature properties of Morris-Thorne wormhole metric}, J. Geom. Phys., \textbf{174(2)} (2022), 104457.
\bibitem{EDS_sultana_2022}
Eyasmin, S., Datta, B. R. and Sarkar, M., \emph{On sultana-dyer spacetime: curvatures and geometric structures}. Int. J. Geom. Methods Mod. Phys., {\bf 20} (2023), 2350101.





\bibitem{ESKiselev2024}
Eyasmin, S., \emph{Curvature properties and Ricci solitons on Kiselev black hole spacetime}. Int. J. Geom. Methods Mod. Phys.,˘(2024), {https://doi.org/10.1142/S0219887825500185}.


\bibitem{F81}
Ferus, D., \emph{A remark on Codazzi tensors on constant curvature space}, Glob. Diff. Geom. Glob. Ann., Lecture
notes 838, Springer, \textbf{1981}.




\bibitem{GG2019}
Gim, Y. and Gwak, B., \emph{Charged particle and strong cosmic censorship in Reissner--Nordstr\"om--de Sitter black holes}, Physical Review D,  amer. phys. soc.,\textbf{100(12)}, 124001 (2019), DOI: 10.1103/PhysRevD.100.124001.



\bibitem{Gray78}
Gray, A., \emph{Einstein-like manifolds which are not Einstein}, Geom. Dedicta, \textbf{7} (1978), 259--280.





\bibitem{Glog02}
G\l ogowska, M.,
\emph{Semi-Riemannian manifolds whose Weyl tensor is a Kulkarni-Nomizu square},
Publ. Inst. Math. (Beograd) (N.S.),
\textbf{72(86)} (2002), 95--106.

\bibitem{Glog-2007}
G\l ogowska, M., \emph{On Roter type manifolds}, Pure and Applied Differential Geometry- PADGE, (2007), 114--122.
	\bibitem{G08}
	G\l ogowska, M., \emph{On quasi-Einstein Cartan type hypersurfaces}, J. Geom. Phys. \textbf{58} (2008), 599--614.

\bibitem{Guler2019}
G\"uler, S. and Crasmareanu, M., \emph{Ricci-Yamabe maps for Riemannian flow and their volume
variation and volume entropy}, Turk. J. Math., \textbf{43}  (2019), 2631--2641.


\bibitem{HV07}
Haesen, S. and Verstraelen, L., \emph{Properties of a scalar curvature invariant depending on two planes}, Manuscripta Math., \textbf{122} (2007), 59-72.
\bibitem{HV07A}
Haesen, S. and Verstraelen, L., \emph{On the scetional curvature of Deszcz}, Anale. Stiint. An. Stiint. Univ. Al. I. Cuza lasi. Mat. (N.S.), \textbf{53} (2007), 181-190.
\bibitem{HV09}
Haesen, S. and Verstraelen, L., \emph{Natural intrinsic geometrical symmetries}, Symmetry, Integrability and Geometry, Methods and Appl. SIGMA, \textbf{5} (2009), 086, 15 pages.
\bibitem{Hamilton1982}
Hamilton, R. S.,
\emph{Three manifolds with positive Ricci curvature},
J. Diff. Geom., 
\textbf{17} (1982), 255--306.


\bibitem{Hamilton1988}
Hamilton, R. S.,
\emph{The Ricci flow on surfaces},
Contemp. Math.,
\textbf{71} (1988), 237--261.

\bibitem{bb7}
 Hawking, S. W. and  Ellis, G. F. R.,  {\it The large scale structure of space-time}, Cambridge University Press, Cambridge \textbf{(1973)}. 
 
\bibitem{Hwkn1975}
Hawking, S. W., \emph{Particle creation by black holes}, Commun.Math. Phys. \textbf{43},199–220 (1975), https://doi.org/10.1007/BF02345020.

\bibitem{bb1}  
Heusler, M., {\it Black hole uniqueness theorems}, Cambridge University Press, Cambridge \textbf{(1996)}. 



%
%





\bibitem{KLD1969}
Katzin, G. H., Livine, J. and Davis, W. R., 
\emph{Curvature collineations: A fundamental symmetry property of the space-times of general relativity defined by the vanishing Lie derivative of the Riemann curvature tensor},
J. Math. Phys., 
\textbf{10(4)}, (1969), 617--629.



\bibitem{KLD1970}
Katzin, G. H., Livine, J. and Davis, W. R.,
\emph{Groups of curvature collineations in Riemannian space-times which admit fields of parallel vectors},
J. Math. Phys.,
\textbf{11} (1970), 1578--1580.


\bibitem{bb6}
 Kerr, R. P., \emph{Gravitational field of a spinning mass as an example of algebraically special metrics}, Phys. Rev. Lett., {\bf 11(5)}, (1963), 237--238. 

%

\bibitem{Kowa06}
Kowalczyk, D., \emph{On the Reissner-Nordström-de Sitter type spacetimes}, Tsukuba J. Math., \textbf{30(2)} (2006), 363--381.



\bibitem{LR89}
Lovelock, D. and Rund, H., \emph{Tensors, differential forms and variational principles}, Courier Dover Publications, \textbf{1989}.

\bibitem{MM12a}
Mantica, C. A. and Molinari, L. G., \emph{Extended Derdzinski-Shen theorem for curvature tensors}, Colloq. Math., \textbf{128} (2012), 1--6.

\bibitem{MM22b}
Mantica, C. A. and Molinari, L. G., \emph{The Jordan algebras of Riemann, Weyl and curvature compatible tensors}, Colloq. Math., \textbf{167} (2022), 63--72.


\bibitem{MM12b}
Mantica, C. A. and Molinari, L. G., \emph{Riemann compatible tensors}, Colloq. Math., \textbf{128} (2012), 197--210.
\bibitem{MM13}
Mantica, C. A. and Molinari, L. G., \emph{Weyl compatible tensors}, Int. J. Geom. Methods Mod. Phys., \textbf{11(08)} (2014), 1450070.

\bibitem{MS12a}
Mantica, C. A. and Suh, Y. J., \emph{The closedness of some generalized curvature 2-forms on a Riemannian manifold I}, Publ. Math. Debrecen, {\bf{81(3-4)}} (2012), 313--326.

\bibitem{MS13a}
Mantica, C. A. and Suh, Y. J., \emph{The closedness of some generalized curvature 2-forms on a Riemannian manifold II}, Publ. Math. Debrecen, {\bf{82(1)}}, (2013), 163--182.
 
\bibitem{MS14}
Mantica, C. A. and Suh, Y. J., \emph{Recurrent conformal 2-forms on pseudo-Riemannian manifolds}, Int. J. Geom. Methods Mod. Phys., \textbf{11(6)} (2014), 1450056 (29 pages).
\bibitem{MS2016}	
Mantica, C. A. and Suh, Y. J., \emph{Pseudo-Z symmetric space-times with divergence-free Weyl tensor and
	pp-waves}, Int. J. Geom. Methods Mod. Phys., \textbf{13(02)} (2016), 1650015.





\bibitem{MIKES76} 
Mike$\check{s}$, J., \emph{Geodesic mappings of symmetric Riemannian spaces}, Odessk. Univ., \textbf{3924-76} (1976), 1-10.
\bibitem{Add1} 
Mike$\check{s}$, J., \emph{Geodesic and holomorphically projective mappings of special Riemannian space}, Ph.D. Thesis,  Odessa Univ.,  (1979).

\bibitem{MS94}
Mike\v s, J. and Sobchuk, V.S., \emph{Geodesic mappings of 3-symmetric Riemannian spaces}
(Russian), translated from Ukrain. Geom. Sb. no. 34 (1991), 80--83, iii J. Math. Sci.
69 (1994), no. 1, 885--887.


\bibitem{MIKES88}
Mike\v s, J., \emph{Geodesic mappings of special Riemannian spaces}, Topics in differential
geometry, Vol. I, II (Debrecen, 1984), Colloq. Math. Soc. J\'anos Bolyai, 46, North-Holland, Amsterdam, (1988), 793--813.


\bibitem{MIKES92}
Mike\v s, J.,\emph{Geodesic mappings of m-symmetric and generalized semisymmetric spaces}
(Russian), translated from Izv. Vyssh. Uchebn. Zaved. Mat. 1992, no. 8, 42-46 Russian Math. (Iz. VUZ), \textbf{36} (1992), no. 8, 38--42.

\bibitem{MIKES96} 
Mike$\check{s}$, J., \emph{Geodesic mapping of affine-connected and Riemannian spaces}, J. Math. Sci., \textbf{78(3)} (1996), 311--333. 

\bibitem{MSV15}
Mike\v s, J., Stepanova, E. and Van\v zurov\'a, A., \emph{ Differential geometry of special mappings},
Palack\'y University Olomouc, Faculty of Science, Olomouc (2015), 568 pp.

\bibitem{MVH09}
Mike\v{s}, J., Van\v{z}urov\'{a}, A. and Hinterleitner, I, \emph{Geodesic mappings and some generalizations}, Palacky Univ. Press, Olomouc, \textbf{2009}.

\bibitem{bb2} 
 Newman, E. T.,  Couch, E., Chinnapared, K.,  Exton, A.,  Prakash, A. and Torrence, R., \emph{Metric of a rotating, charged mass}, J. Math. Phys., {\bf (6)}, (1965), 918--919. 



\bibitem{Nord1918}
Nordström, G., \emph{On the Energy of the gravitational field in Einstein’s theory}, Verhandl. Koninkl. Ned. Akad. Wetenschap., Afdel. Natuurk., \textbf{26} (1918), 1201-1208.

\bibitem{bb5} Nordstr\"{o}m, G., \emph{On the energy of the
gravitational field in Einstein's theory}, Proc. Kon. Ned. Akad. Wet., {\bf 20(2)} (1918), 1238--1245. 






\bibitem{VBdS_1987}
  Patino, A. and Rago, H., \emph{A radiating charge embedded in a De Sitter universe},   Phys. Lett. A, \textbf{121(7)} (1987), 329--330.
  
\bibitem{Neill83}
O'Neill, B., \emph{Semi-Riemannian geometry with applications to the relativity}, Academic Press, New York-London, \textbf{1983}.
\bibitem{Patt52} 
Patterson, E. M., \emph{Some theorems on Ricci recurrent spaces}, J. London Math. Soc., \textbf{27} (1952), 287--295.
\bibitem{Pigola2011}
Pigola, S., Rigoli, M., Rimoldi, M., Setti, A. G., 
\emph{Ricci almost solitons},
Ann. Scuola Norm. Sup. Pisa Cl. Sci., \textbf{X(5)} (2011), 757--799.




\bibitem{P95}
Prvanovi$\acute{\mbox{c}}$,  M.,  \emph{On weakly symmetric Riemannian manifolds,} Publ. Math. Debrecen, \textbf{46(1-2)} (1995), 19--25.



\bibitem{Reis1916}
Reissner, H., \emph{Über die Eigengravitation des elektrischen Feldes
	nach der Einsteinschen Theorie}, Anna. der Physik, \textbf{50} (1916), 106--120.

\bibitem{Ruse46}
Ruse, H. S., \emph{On simply harmonic spaces}, J. London Math. Soc., \textbf{21} (1946), 243--247.
\bibitem{Ruse49a}
Ruse, H. S., \emph{On simply harmonic `kappa spaces' of four dimensions}, Proc. London Math. Soc., \textbf{50}
(1949), 317--329.
\bibitem{Ruse49b}
Ruse, H. S., \emph{Three dimensional spaces of recurrent curvature}, Proc. London Math. Soc., \textbf{50} (1949),
438--446.





\bibitem{Sch1916}
Schwarzschild, K., \emph{Über das Gravitationsfeld eines Massenpunktes nach der Einsteinschen Theorie}, Sitzungsber. Preuss. Akad. D. Wiss. \textbf{50} (1916), 189--196.
\bibitem{S09}
Shaikh, A. A., \emph{On pseudo quasi-Einstein manifolds}, Period. Math. Hungarica, \textbf{59(2)} (2009), 119--146.

\bibitem{SAA18}
Shaikh, A. A., Ali, M. and Ahsan, Z., \emph{Curvature properties of Robinson-Trautman metric}, J. Geom., \textbf{109(2)} (2018), 1--20. DOI: 10.1007/s00022-018-0443-1



\bibitem{SAAC20}
Shaikh, A. A., Ali, A., Alkhaldi, A. H. and Chakraborty, D., \emph{Curvature properties of Melvin magnetic metric}, J. Geom. Phys., \textbf{150} (2020), 103593. DOI: 10.1016/j.geomphys.2019.103593.

\bibitem{SAAC20N}
Shaikh, A. A., Ali, A., Alkhaldi, A. H. and Chakraborty, D., \emph{Curvature properties of Nariai spacetimes}, Int. J. Geom. Methods Mod. Phys., \textbf{17(03)} (2020), 2050034. DOI: 10.1142/S0219887820500346

\bibitem{SAACD_LTB_2022}
Shaikh, A. A., Ali, A., Alkhaldi, A. H., Chakraborty, D. and Datta, B. R., \emph{On some curvature properties of Lemaitre–Tolman–Bondi spacetime}. Gen. Relativ. Gravit. \textbf{54(1)} (2022), 6 (21 pages). DOI: 10.1007/s10714-021-02890-4

\bibitem{SAR13}
Shaikh, A. A., Al-Solamy, F. R. and Roy, I., 
\emph{On the existence of a new class of semi-Riemannian manifolds}, 
Mathematical Sciences, \textbf{7} (2013), 46.



\bibitem{SASZ2022} 
Shaikh, A. A., Ali, M., Salman, M. and Zengin, F. O., \emph{Curvature inheritance symmetry on M-projectively flat spacetimes}, Int. J. Geom. Methods Mod. Phys., \textbf{20(2)}  (2023), 2350088.


\bibitem{SADM_TCEBW_2023}
Shaikh, A. A., Ahmed, F., Datta, B. R. and Sarkar, M.,
\emph{On geometric groperties of topologically charged Ellis-Bronnikov-type wormhole spacetime}, Filomat, {\bf 38(15)} (2024). 

\bibitem{SAS_EiBI_2024}
Shaikh, A. A., Ahmed, F. and Sarkar, M., 
\emph{Curvature related geometrical properties of topologically charged EiBI-gravity spacetime}, New Astronomy, {\bf 112} (2024), 102272 . 

\bibitem{SASK_pgm_2024}
Shaikh, A. A., Ahmed, F., Sarkar, M. and Kamiruzzaman, 
\emph{Symmetry and pseudosymmetry properties of a point-like global monopole spacetime
}, Int. J. Geom. Methods Mod. Phys.,
 \textbf{(2025)}, 2550156, DOI:10.1142/S0219887825501567.	

\bibitem{SASK_LBH_2025}
Shaikh, A. A., Ahmed, F., Sarkar, M. and Kamiruzzaman, 
\emph{An exploration of geometric and curvature properties of Lemos black hole spacetime
}, Chinese J. Phys.,
 \textbf{96}, (2025), 333--354.

\bibitem{SB08}
Shaikh, A. A. and Binh, T. Q., \emph{On some class of Riemannian manifolds,} Bull. Transilv. Univ., \textbf{15(50)} (2008), 351--362.
\bibitem{SBH21}
Shaikh, A. A., Binh, T. Q. and Kundu, H., \emph{Curvature properties of generalized $pp$-wave metrics}, Kragujevac J. Math.,
\textbf{45(2)} (2021), 237--258.


\bibitem{SC21}
Shaikh, A. A. and Chakraborty, D., \emph{Curvature properties of Kantowski-Sachs metric}, J. Geom. Phys., \textbf{160} (2021), 103970. DOI: 10.1016/j.geomphys.2020.103970


\bibitem{SDAA_LCS_2021}
Shaikh, A. A., Datta, B. R., Ali, A. and Alkhaldi, A. H., \emph{LCS-manifolds and Ricci solitons}. Int. J. Geom. Methods Mod. Phys., \textbf{18(09)} (2021), 2150138.

\bibitem{ShaikhDatta2022} Shaikh, A. A. and Datta, B. R., \emph{Ricci solitons and curvature inheritance on Robinson-Trautman spacetimes}, Int. J. Geom. Methods Mod. Phys., \textbf{21(99)} (2024), https://doi.org/10.1142/S0219887824501639.



\bibitem{SDC}
Shaikh, A. A., Datta, B. R. and Chakraborty, D., \emph{On some curvature properties of Vaidya-Bonner metric}, Int. J. Geom. Methods. Phys., https://doi.org/10.1142/S0219887821502054.


\bibitem{SDKC19}
Shaikh, A. A., Das, L., Kundu, H. and Chakraborty, D., \emph{Curvature properties of Siklos metric}, Diff.
Goem.- Dyn. Syst., \textbf{21} (2019), 167--180.




\bibitem{SDHJK15}
Shaikh, A. A., Deszcz, R., Hotlo\'{s}, M., Je\l owicki, J. and Kundu, H.,
\emph{On pseudosymmetric manifolds},
Publ. Math. Debrecen,
\textbf{86(3-4)} (2015), 433-456.


\bibitem{SHDS_hayward}
Shaikh, A. A., Hui, S. K., Datta, B. R. and Sarkar, M.,
\emph{On curvature related geometric properties of Hayward black hole spacetime}, New Astronomy, {\bf 108} (2024), 102181. 
	 	

\bibitem{SHS_warped}
Shaikh, A. A., Hui, S. K. and Sarkar, M., \emph{Curvature properties of a warped product metric},  Palestine J. Math., \textbf{13(1)} (2024), 220--231.
		
\bibitem{SHS_Bardeen}
Shaikh, A. A., Hui, S. K. and Sarkar, M., \emph{Curvature properties of Bardeen black hole spacetime}, Bulgarian J. Phys., \textbf{50} (2023), 168. 



\bibitem{SHSB_VBDS}
Shaikh, A. A., Hui, S. K., Sarkar, M. and Babu, V. A.,  \emph{Symmetry and pseudosymmetry properties of Vaidya-Bonner-de Sitter spacetime}, J. Geom. Phys., \textbf{202} (2024), 105235. 



\bibitem{SJ06}
Shaikh, A. A. and Jana, S. K., \emph{On weakly cyclic Ricci symmetric manifolds,} Ann. Pol. Math., \textbf{89(3)} (2006), 139--146.

\bibitem{SJ07}
Shaikh, A. A. and Jana, S. K., \emph{On quasi-conformally flat weakly Ricci symmetric manifolds,} Acta Math. Hungar., \textbf{115(3)} (2007), 197--214.

\bibitem{SKH11}
Shaikh, A. A., Kim, Y. H. and Hui, S. K., \emph{On Lorentzian quasi Einstein manifolds}, J. Korean Math. Soc., \textbf{48} (2011), 669--689. Erratum: \emph{On Lorentzian quasi Einstein manifolds,} J. Korean Math. Soc., \textbf{48(6)} (2011), 1327--1328.


\bibitem{SK14}
Shaikh, A. A. and Kundu, H., \emph{On equivalency of various geometric structures}, J. Geom., \textbf{105} (2014), 139--165. DOI: 10.1007/s00022-013-0200-4

\bibitem{SK16}
Shaikh, A. A. and Kundu, H., \emph{On warped product generalized Roter type manifolds}, Balkan J. Geom. Appl., \textbf{21(2)} (2016), 82--95.

\bibitem{SK16srs}
Shaikh, A. A. and Kundu, H., \emph{On curvature properties of Som-Raychaudhuri spacetime}, J. Geom., \textbf{108(2)} (2016), 501--515.
\bibitem{SKppsnw}
Shaikh, A. A. and Kundu, H., \emph{On warped product manifolds satisfying some pseudosymmetric type conditions}, Diff. Geom. - Dyn. Syst., \textbf{19} (2017), 119--135.


\bibitem{SK19}
Shaikh, A. A. and Kundu, H., \emph{On generalized Roter type manifolds},
Kragujevac J. Math., \textbf{43(3)} (2019), 471--493.


\bibitem{SKA18}
Shaikh, A. A., Kundu, H. and Ali, Md. S., \emph{On warped product super generalized recurrent manifolds}, An. \c{S}tiin\c{t}. Univ. Al. I. Cuza Ia\c{s}i. Mat. (N. S.), \textbf{LXIV(1)} (2018), 85--99.

\bibitem{SKS19}
Shaikh, A. A., Kundu, H. and Sen, J., \emph{Curvature properties of Vaidya metric}, Indian J. Math., \textbf{61(1)} (2019), 41--59.
\bibitem{SMM2022} Shaikh, A. A., Mandal, P. and  Mondal, C. K., {\em Diameter estimation of gradient $\rho$-Einstein soliton}, J. Geom. Phys., \textbf{177} (2022), 104518.

\bibitem{SM23}Shaikh, A. A. and Mandal, P., \textit{Some characterizations of quasi Yamabe solitons}, J. Geom., \textbf{114} (2023).
\bibitem{Absos23}Shaikh, A. A., Mandal, P. and Mondal, C. K., {\em Existence of finite global norm of potential vector field in a Ricci Soliton}, Proc. Natl. Acad. Sci., India, Sect. A Phys. Sci., \textbf{93} (2023), 671--674.

\bibitem{SMB24} Shaikh, A. A., Mandal, P., Babu, V. A., \textit{Triviality results and conjugate radius estimation of Ricci solitons}, Bull. Braz. Math. Soc. (N. Ser.), \textbf{55(22)} (2024).
 
\bibitem{SP10}
Shaikh, A. A. and Patra, A., \emph{On a generalized class of recurrent manifolds}, 
Arch. Math. (BRNO), {\bf 46} (2010), 71--78.
\bibitem{SR11}
Shaikh, A. A. and Roy, I., \emph{On weakly generalized recurrent manifolds}, Ann. Univ. Sci. Budapest, E$\ddot{\mbox{o}}$tv$\ddot{\mbox{o}}$s Sect. Math., \textbf{54} (2011), 35--45.
\bibitem{SRK18}
Shaikh, A. A., Roy, I. and Kundu, H.,
\emph{On the existence of a generalized class of recurrent manifolds},  An. \c{S}tiin\c{t}. Univ. Al. I. Cuza Ia\c{s}i. Mat. (N. S.), \textbf{LXIV(2)} (2018), 233--251.
\bibitem{SRK15}
Shaikh, A. A., Roy, I. and Kundu, H.,
\emph{On the existence of a generalized class of recurrent manifolds},  An. \c{S}tiin\c{t}. Univ. Al. I. Cuza Ia\c{s}i. Mat. (N. S.), \textbf{LXIV(2)} (2018), 233--251.

\bibitem{SRK17}
Shaikh, A. A., Roy, I. and Kundu, H., \emph{On some generalized recurrent manifolds}, Bull.
Iranian Math. Soc., \textbf{43(5)} (2017), 1209--1225.



\bibitem{SSC19}
Shaikh, A. A., Srivastava, S. K. and Chakraborty, D., \emph{Curvature properties of anisotropic scale invariant
metrics}, Int. J. Geom. Meth. Mod. Phys., {\bf 16} (2019), 1950086.

\bibitem{SYH09}
Shaikh, A. A., Yoon, D. W. and Hui, S. K., \emph{ On quasi-Einstein spacetimes}, Tsukuba J. Math., {\bf 33(2)} (2009), 305--326.
\bibitem{Shi25}
Shirokov, P. A., \emph{Constant fields of vectors and tensors of second order on Riemannian spaces}, Kazan, U\v cen. Zap. Univ.,  {\bf 25(2)} (1925), 256--280.


\bibitem{Shi98}
Shirokov, P. A., \emph{Shirokov's work on the geometry of symmetric spaces}, J. Math. Sci., (New York), {\bf 89(3)} (1998), 1253--1260.

\bibitem{Siddiqi2020}
Siddiqi, M. and Akyol, M. A.,  \emph{$\eta $-Ricci-Yamabe soliton on Riemannian submersions from Riemannian manifolds}. arXiv preprint arXiv:2004.14124, (2020).

\bibitem{S81}
Simon, U., \emph{Codazzi tensors}, Glob. Diff. Geom. and Glob. Ann., Lecture notes, 838, Springer-Verlag, 1981,
289--296.
\bibitem{Sin54}
Sinyukov, N. S., \emph{On geodesic mappings  of Riemannian manifolds onto symmetric spaces}, Dokl. Akad. Nauk SSSR,  {\bf 98} (1954),
21--23.
\bibitem{SG1986}	
Sippel, R. and Goenner, H., \emph{Symmetry classes of pp-waves}, Gen. Relativ. Gravit. \textbf{18(12)} (1986), 1229--1243.

\bibitem{SKMHH03}
Sthepani, H., Kramer, D., Mac-Callum, M., Hoenselaers, C. and Hertl, E., \emph{Exact solutions
of Einstein's field equations}, Cambridge Monographs on Mathematical Physics, Cambridge
University Press, Second Edition, \textbf{2003}.
\bibitem{SKP03}
Suh, Y. J., Kwon, J-H. and Pyo, Y. S., \emph{On semi-Riemannian manifolds satisfying the second Bianchi identity}, J. Korean Math. Soc., \textbf{40(1)} (2003), 129--167.
\bibitem{Szab82}
Szab$\acute{\mbox{o}}$, Z. I., \emph{Structure theorems on Riemannian spaces satisfying $R(X,Y)\cdot R=0$, I. The local version}, J. Diff. Geom., \textbf{17} (1982), 531--582.

\bibitem{Szab84}
Szab\'o, Z. I., \emph{Classification and construction of complete hypersurfaces satisfying 
$R(X, Y)\cdot R = 0$}, Acta Sci. Math., \textbf{47} (1984), 321--348.

\bibitem{Szab85}
Szab\'o, Z. I., \emph{Structure theorems 
on Riemannian spaces satisfying $R(X, Y)\cdot R = 0$, II, The global version}, 
Geom. Dedicata, \textbf{19} (1985), 65--108.



\bibitem{Tach74}
Tachibana, S.,
\emph{A theorem on Riemannian manifolds of positive curvature operator},
Proc. Japan Acad.,
\textbf{50} (1974), 301--302.






\bibitem{Venz85} 
Venzi, P., 
\emph{Una generalizzazione degli spazi ricorrenti}, 
Rev. Roumaine Math. Pures Appl.,
\textbf{30} (1985), 295--305.

\bibitem{Walk50}
Walker, A. G.,
\emph{On Ruse's spaces of recurrent curvature},
Proc. London Math. Soc.,
\textbf{52} (1950), 36--64.
\bibitem{Zhao2014}
Zhao, et. al., \emph{The critical phenomena and thermodynamics of the Reissner-Nordstr\"om-de Sitter black hole}, Adv. High Energy Phys, \textbf{2014}, 124854 (2014).

\bibitem{Zhen2025}
Zhen et. al., \emph{Schottky anomaly of Reissner-Nordstr\"om-de Sitter
spacetime}, Chinese Phys. C \textbf{49(3)}, (2025) 035105.
{\color{blue}









 














































}



\end{thebibliography}
\end{document}